\documentclass[11pt]{article}

\usepackage{amssymb,latexsym,amsfonts,verbatim,amscd}

\newtheorem{Def}{Definition}[section]
\newtheorem{Thm}[Def]{Theorem}
\newtheorem{Lem}[Def]{Lemma}
\newtheorem{Prop}[Def]{Proposition}
\newtheorem{Cor}[Def]{Corollary}
\newtheorem{Rem}[Def]{Remark}
\newtheorem{Rems}[Def]{Remarks}
\newtheorem{Fac}[Def]{Fact}

\newtheorem{Ex}[Def]{Example}
\newtheorem{Quest}[Def]{Question}
\font\nat msbm10 scaled\magstephalf
\def\N{\hbox{\nat\char78}}

\font\mata=msam10 
\def\restr{\mbox{\mata\char22}}
\def\telos{\hfill$\dashv$}

\begin{document}

\title{Propositional superposition logic}
\author{Athanassios Tzouvaras}

\date{}
\maketitle

\begin{center}
Department  of Mathematics\\
Aristotle University of Thessaloniki \\
541 24 Thessaloniki, Greece \\
e-mail: \verb"tzouvara@math.auth.gr"
\end{center}

\begin{abstract} We extend classical Propositional Logic (PL) by adding a new primitive binary  connective $\varphi|\psi$, intended to represent the ``superposition'' of sentences $\varphi$ and $\psi$, an operation motivated  by the corresponding notion of  quantum mechanics, but not intended to  capture all aspects of the latter as they appear in physics. To interpret the new connective, we extend the classical Boolean semantics by employing models  of the form $\langle M,f\rangle$, where $M$ is an ordinary  two-valued assignment for the sentences of PL and $f$ is  a choice function  for all pairs of classical sentences. In the new semantics $\varphi|\psi$ is  strictly interpolated between $\varphi\wedge\psi$ and $\varphi\vee\psi$. By imposing several constraints on the choice functions we obtain corresponding notions of  logical consequence  relations and corresponding systems of tautologies, with respect to which $|$ satisfies some natural algebraic properties such as associativity, closedness under logical equivalence and distributivity over its dual connective. Thus various systems of Propositional Superposition Logic (PLS) arise as extensions of PL. Axiomatizations for these systems of tautologies are presented and soundness is shown  for all of them. Completeness is proved for the weakest of these systems. For the other  systems  completeness holds if and only if  every consistent set of sentences is extendible to a consistent and complete one, a condition whose truth is closely related to the validity of the deduction theorem.
\end{abstract}

\vskip 0.2in


\vskip 0.2in

 Choice function for pairs of sentences.  Associative, regular,  $\neg$-decreasing choice functions/orderings.  Propositional superposition logics (PLS).

\newpage

\section{Introduction}
In this paper we present an extension  of classical Propositional Logic (PL) (more precisely, an array of extensions of increasing strength), obtained by adding a new binary logical operation, called ``superposition'', together with a new semantics extending the standard one, inspired and motivated by the corresponding notion of  quantum mechanics.  That the notion of superposition is central  in quantum mechanics is rather well-known. For the sake of completeness let us outline briefly the core of the  idea.

A quantum system $A$,  for example   an electron or  a photon,  can  be only in finitely many  possible ``states'' (or rather ``pure states'', see \cite{DaGu02})   with respect to a certain physical magnitude $Q$ such as spin,  charge, etc. Suppose for simplicity  that $A$ can be only in two  possible states, say ``spin up'' and ``spin down''.  We know that whenever the spin of $A$ is measured, the outcome will necessarily be   either ``spin up'' or ``spin down'' but one cannot predict it in advance  precisely, except only with  a certain degree of probability. While unobserved, $A$ is thought to be at some kind of a  mixture or composition of  these states,   called {\em  superposition of states.} But as soon as the system $A$ is scanned  and measured, the superposition breaks down, or ``collapses''  according to the jargon of quantum mechanics, to one of the constituent states. So in a sense,   the states ``spin up'' and ``spin down''  co-exist and at the same time  exclude  each other.

From its very beginning quantum mechanics had developed an effective and flexible  formalism  to represent the states of a system (see \cite{DaGu02} for a brief overview of the subject), namely as  vectors of a Hilbert space. For example the pure states  ``spin up'' and ``spin down'' are represented by vectors $\vec{u}_0$,  $\vec{u}_1$, respectively. Then the ``principle of superposition'' says that for any complex numbers $c_0$, $c_1$ such that $|c_0|^2+|c_1|^2=1$, the linear  combination  $c_0\vec{u}_0+c_1\vec{u}_1$ is also a legitimate state of the system $A$. Moreover,  $|c_0|^2$, $|c_1|^2$ represent the probabilities for  $A$ to be in  state $\vec{u}_0$ or   $\vec{u}_1$, respectively, when measured.
This  treatment of superposition as a linear  combination of vectors is mainly due to P.A.M.  Dirac\footnote{P.A.M. Dirac, {\em The Principles of Quantum Mechanics,} Oxford U.P., 1958.}, who considered the principle of superposition as one of the most fundamental properties of quantum mechanics.

Later on a new  approach to quantum mechanics through {\em quantum logics} was developed  by  the  work of G. Birkhoff- J. von Neumann\footnote{ G. Birkhoff and J. von Neumann, The logic of quantum mechanics, {\em Ann. Math.} {\bf 37} (1936), 823-843.}, G. Mackey\footnote{G. Mackey, {\em The Mathematical Foundations of Quantum Mechanics,} W. A. Benjamin, 1963.} and others. Here the emphasis was in the  formalization  of {\em non-distributivity,} another characteristic phenomenon of quantum mechanics, and it was not clear  whether and how non-distributivity and superposition were related to each other.  As S. Gudder says\footnote{S.P. Gudder, A superposition principle in physics, {\em J. Math. Phys.} {\bf 11} (1970), no. 3, 1037-1040.},  the problem arose to find a formulation of the principle of superposition in the quantum logic  approach, roughly equivalent to Dirac's formulation in the vector-space approach.  Various such formulations of superposition can be found in the literature.\footnote{See for example S.P.Gudder above  and S. Pulmannov\'{a}, A superposition principle in quantum logics, {\em Commun. Math. Phys.} {\bf 49} (1976), no. 3, 47-51.}  Note that all versions of quantum logic are weaker than classical logic, since they lack the distributivity law.

An important source of inspiration for the present work has been   E. Schr\"{o}dinger's  1935 paper \cite{Sch35} containing the ``cat paradox'', in which the author shows, by his famous  thought experiment, how   superposition of quantum states might (in principle) be transformed into superposition of {\em macroscopic situations.} Although Schr\"{o}dinger refers to the experiment  with ridicule, as a ``serious misgiving arising if one notices that the uncertainty affects macroscopically tangible and visible things'', it is perhaps the first hint towards thinking that the  phenomenon could be conceived in a much broader sense, even in contexts different from the original one. And it is this general,  abstract and purely {\em logical content} of superposition that we are interested in and  deal with in this paper.\footnote{As already said earlier, whether the logical content of superposition, as isolated here, bears  actual {\em connections} with and/or {\em applications} to the existing systems of quantum mechanics and quantum logic is not known at present. Some further comments on this issue are made in the last section.}

In particular, the purpose of the paper is  to offer a simple  interpretation of superposition not by means of a variant of quantum logic, but rather by an  {\em extension of classical logic.}  The interpretation is absolutely within  classical reasoning  and common sense,  since we do not drop any law of classical logic, but only augment them by new ones concerning the  superposition operation. The  ingredient that makes it possible to go  beyond classical tautologies is the use at each truth evaluation of a {\em choice function} acting upon pairs of sentences, a tool originating in set-theory rather than logic.

Let $\varphi_0$, $\varphi_1$ denote the statements ``$A$ is at state $\vec{u}_0$'' and ``$A$ is at state $\vec{u}_1$'', respectively, and let $\varphi_0|\varphi_1$ denote the statement  ``$A$ is at the  superposition of states $\vec{u}_0$ and  $\vec{u}_1$''.  $\varphi_0$, $\varphi_1$ are ordinary  statements, so they can be assigned ordinary  truth values. But what about the truth values of $\varphi_0|\varphi_1$? Clearly the operation  $\varphi_0|\varphi_1$ cannot be  expressed  in classical logic, that is,  $\varphi_0|\varphi_1$ cannot be logically equivalent to a Boolean combination  $S(\varphi_0,\varphi_1)$ of $\varphi_0$, $\varphi_1$.\footnote{As is well-known there exist precisely 16  classical binary operations $S(\varphi_0,\varphi_1)$, definable in terms of $\wedge$, $\vee$ and $\neg$ (including  $\wedge$, $\vee$ themselves, and also  $\rightarrow$, $\leftrightarrow$, their negations, as well as other trivial ones), none of which can express the logical content of $|$. } However, an  intriguing
feature of $\varphi_0|\varphi_1$  is that it has points in common  {\em both} with classical conjunction and classical disjunction. In a sense it is a ``mixture'' of  $\varphi_0\wedge \varphi_1$ and $\varphi_0\vee \varphi_1$, or a property between them, since it bears a  conjunctive as well  as a disjunctive component. Indeed, $\varphi_0|\varphi_1$ means on the one hand  that the properties $\varphi_0$ and $\varphi_1$ hold {\em simultaneously} (at least partly) during the non-measurement phase, which is clearly  a conjunctive component of $\varphi_0|\varphi_1$, and on the other, at any particular collapse of the superposed states during a measurement, $\varphi_0|\varphi_1$ reduces to either $\varphi_0$ {\em or} $\varphi_1$, which is  a disjunctive component of the operation. The interpretation of $\varphi_0|\varphi_1$ given below  justifies in fact this meaning of $\varphi_0|\varphi_1$ as ``something between  $\varphi_0\wedge \varphi_1$ and $\varphi_0\vee \varphi_1$''.

Let us consider  a propositional language $L=\{p_0,p_1,\ldots\}\cup\{\wedge,\vee,\neg\}$, where $p_i$ are symbols of atomic propositions, whose interpretations are  usual  two-valued truth assignments $M:Sen(L)\rightarrow \{0,1\}$ respecting the classical truth tables. Let us extend $L$ to $L_s=L\cup\{|\}$, where $|$ is a new primitive binary connective. For any sentences $\varphi,\psi$ of $L_s$, $\varphi|\psi$  denotes the superposition of $\varphi$ and $\psi$.
Then an interpretation for the sentences of $L_s$ can be given by the help of  a truth assignment  $M$ for the sentences of $L$, together with a {\em collapsing   mapping}  $\textsf{c}$ from the sentences  of $L_s$ into those of $L$. The mapping $\textsf{c}$ is intended to represent the collapsing of the superposed $\varphi|\psi$ to one of its components.

The basic idea is that the collapse  of the composite state  $c_0\vec{u}_0+c_1\vec{u}_1$ to one of the sates $\vec{u}_0$, $\vec{u}_1$  can be seen, {\em from the point of view of pure logic,} just as a (more or less random)  choice from  the set of possible outcomes $\{\vec{u}_0,\vec{u}_1\}$. This is because from the point of view of pure logic   probabilities are  irrelevant or, which amounts to the same thing, the states $\vec{u}_0$ and $\vec{u}_1$ are considered  equiprobable.   In such a case  the  superposition of $\vec{u}_0$ and $\vec{u}_1$ is unique and the outcome of the collapse can be decided  by a coin tossing or, more strictly,  by  a  {\em choice function} acting on  pairs of observable states, which in our case coincide with pairs of sentences of  $L$. This  of course constitutes  a major deviation from the standard treatment of superposition, according to which there is not just one superposition of  $\vec{u}_0$ and  $\vec{u}_1$ but infinitely many, actually as many  as the number of  linear combinations $c_0\vec{u}_0+c_1\vec{u}_1$, for $|c_0|^2+|c_1|^2=1$. So  the logic presented here is hardly  the logic of superposition as this concept is currently used and understood in  physics today. It is rather the logic of superposition, when the latter is  understood as the ``logical extract'' of the corresponding physics concept. Whether it could eventually have applications to the field of quantum mechanics we don't know.

The elementary requirements for a collapsing map  $\textsf{c}$ are the following: (a) it must be the identity on classical sentences, that is, $\textsf{c}(\varphi)=\varphi$ for every   $L$-sentence $\varphi$. (b) It  must commute with the standard connectives $\wedge$, $\vee$ and $\neg$,  that is, $\textsf{c}(\varphi\wedge \psi)=\textsf{c}(\varphi)\wedge \textsf{c}(\psi)$, $\textsf{c}(\varphi\vee \psi)=\textsf{c}(\varphi)\vee \textsf{c}(\psi)$ and $\textsf{c}(\neg\varphi)=\neg\textsf{c}(\varphi)$. (c)  $\textsf{c}(\varphi|\psi)$ must be   {\em some} of the sentences $\textsf{c}(\varphi)$, $\textsf{c}(\psi)$, which is chosen by the help of a choice function $f$ for pairs of classical sentences, that is, $$\textsf{c}(\varphi|\psi)=f(\{\textsf{c}(\varphi),\textsf{c}(\psi)\}).$$
Since every sentence of $L_s$ is built from atomic sentences all of which belong to the initial classical language $L$, it follows that $\textsf{c}$ is {\em fully determined} by the choice function $f$, and below we shall write  $\textsf{c}=\overline{f}$. Therefore choice functions $f$ for  pairs of sentences of $L$  are the cornerstone of the new semantics.

Given a truth assignment $M$ for $L$ and a choice function  $f$ for $L$, a sentence  $\varphi$ of $L_s$ is {\em true in $M$ under the  choice function  $f$,} denoted $\langle M,f\rangle\models_s \varphi$, if and only if  $\overline{f}(\varphi)$ is (classically) true in $M$. That is:
$$\langle M,f\rangle\models_s \varphi \  \mbox{iff} \ M\models \  \overline{f}(\varphi).$$
Since $\overline{f}$ is generated by $f$, special conditions on $f$ induce special properties for $\overline{f}$ that in turn affect the properties of $\models_s$. Such a condition is needed, for instance, in order for $|$ to be associative.

The above truth concept $\models_s$   extends the classical one and  induces the  notions of s-logical consequence, $\varphi\models_s\psi$, and s-logical equivalence $\sim_s$,  which generalize the corresponding standard  relations $\models$ and $\sim$. A nice feature of the new semantics is that for all sentences  $\varphi$, $\psi$,
\begin{equation} \label{E:between}
\varphi\wedge \psi \models_s  \varphi|\psi \models_s  \varphi\vee \psi,
\end{equation}
where the relations $\models_s$ in both places are strict, that is, they cannot in general be reversed (see Theorem \ref{T:interpol} below).  It means that  $\varphi|\psi$ is {\em strictly interpolated} between $\varphi\wedge \psi$ and $\varphi\vee \psi$, a fact that in some sense makes precise the above expressed intuition that $\varphi|\psi$ is a ``mixture'' of $\varphi\wedge \psi$ and $\varphi\vee \psi$. In particular, $$\varphi\wedge \neg\varphi \models_s  \varphi|\neg\varphi \models_s  \varphi\vee \neg\varphi,$$
which means that the superposition of two contradictory situations, like those in Schr\"{o}dinger's cat experiment \cite{Sch35} mentioned above,  is neither  a contradiction nor a paradox at all (see Corollary \ref{C:Sch} below).

Another nice feature of the semantics is that in order for $|$ to be associative with respect to a structure $\langle M,f\rangle$, that is, $\langle M,f\rangle\models_s\varphi|(\psi|\sigma)\leftrightarrow (\varphi|\psi)|\sigma$, it is necessary and sufficient for $f$ to coincide with the function $\min_<$ induced by a total ordering $<$ of the set of sentences of $L$. Such an  $f=\min_<$ picks from each pair $\{\alpha,\beta\}$ not a ``random''  element  but  the {\em least} one with respect to  $<$. This kind of choice functions  will be  the
dominant one throughout the paper.

No knowledge of quantum mechanics or quantum logic is required for reading this paper. The only prerequisite  is just  knowledge of basic  Propositional Logic (PL), namely its semantics, axiomatization and soundness and completeness theorems, as well as some elementary set-theoretic facts concerning choice functions for sets of finite sets, total orderings etc.  For example \cite{En02} is one of the many logic texts that contain  the necessary  material.  Nevertheless, some familiarity with non-classical logics, their axiomatization  and their semantics,  would be highly  helpful. Also for the subject of  choice functions and choice principles the reader may consult \cite{Je73}.

Finally, I should mention some  other current treatments of superposition from a logical point of view, although one can hardly find to them points of overlapping and convergence  with the present one. Such logical approaches are contained in \cite{CoRo13}, \cite{KrAr15} and \cite{Ba10},  to mention the most recent ones. The main difference of theses approaches  from the present one is that they are all based on some non-classical logical system, while our point of departure is the solid ground of classical propositional logic. For instance \cite{CoRo13} relies heavily on paraconsistent logic that allows one to  accommodate contradictions without collapsing the system. In fact  superposition is captured in  \cite{CoRo13} as a ``contradictory situation'': if a quantum system $S$ is in the state of superposition of the states $s_1$ and $s_2$, this is expressed by the help of a two-place predicate $K$ and the conjunction of  axioms $K(S,s_1)$, $\neg K(S,s_1)$, $K(S,s_2)$ and  $\neg K(S,s_2)$. (Here the negation $\neg$ is ``weak'' and the conjunction of these claims is not catastrophic.) Analogously,  \cite{KrAr15} uses a version of modal logic in an enriched language that, besides $\neg$, $\wedge$, $\vee$ and $\diamond$ (possibility operator), contains a binary connective $\star$ for the superposition operation and a unary connective $M$ for ``measurement has been done''. Also a  Kripke semantics is used, and the basic idea, as I understood it, is  to avoid the contradiction arising e.g.  from Schr\"{o}dinger's cat, by ``splitting'' it, after the measurement, between two different possible worlds, one containing the cat alive and one containing the cat dead.  Finally  \cite{Ba10} is more syntactically oriented. It treats superposition syntactically by employing a version of sequent calculus called ``basic logic'' (developed in \cite{SBF00}), which encompasses aspects of linear logic and quantum logic.

\vskip 0.1in

{\bf Summary of Contents.} Section 2 contains the semantics of $|$ based on choice functions for pairs of sentences of $L$. More specifically, in subsection 2.1 we give the basic definitions of the new semantics and the corresponding notions of logical consequence $\models_s$ and  logical equivalence $\sim_s$. The models for the sentences of $L_s$ are structures of the form $\langle M,f\rangle$, where $M$ is a truth assignment to sentences of $L$ and $f$ is an arbitrary choice function for $L$.
We prove the basic facts, among which that $\varphi\wedge\psi\models_s\varphi|\psi\models_s\varphi\vee\psi\models_s$.
The properties of $|$ supported by such general structures are only $\varphi|\varphi\leftrightarrow\varphi$ (idempotence) and $\varphi|\psi\leftrightarrow\psi|\varphi$ (commutativity). In order to obtain further properties for $|$ we need to impose additional conditions on the choice functions employed which entail more and more refined truth  notions.  In general if ${\cal F}$ is the set of all choice functions for $L$, for any nonempty $X\subseteq {\cal F}$  the relations $\models_X$, of $X$-logical consequence, and $\sim_X$, of $X$-logical equivalence,  are naturally defined by employing  models $\langle M,f\rangle$ with $f\in X$ (rather than $f\in {\cal F}$). For each such $X\subseteq {\cal F}$ the set of $X$-tautologies $Taut(X)=\{\varphi:\models_X\varphi\}$ is defined.

In  the next subsections of \S 2 we focus on certain natural such subclasses $X\subseteq {\cal F}$ and the corresponding truth notions.

In subsection 2.2 we introduce the class $Asso$ of {\em associative} choice functions and a  simple and elegant characterization of them is given, as the functions $\min_<$ with respect to total orderings $<$ of the set of $L$-sentences. The term comes from the fact that if $f\in Asso$, then  $|$ is associative with respect to every structure  $\langle M,f\rangle$.  A kind of converse holds also: If $|$ is associative with respect to $\langle M,f\rangle$, then $f$ is ``essentially associative''.

In subsection 2.3 we introduce the class $Reg$ of {\em regular} choice functions,  as well as the finer  class $Reg^*=Reg\cap Asso$. Regularity guarantees that the truth relation $\models_{Reg}$, as well as $\models_{Reg^*}$, is ``logically closed'', that is, for any subsentence $\sigma$  of $\varphi$ and any $\sigma'\sim_{Reg}\sigma$, $\varphi[\sigma'/\sigma]$ and $\varphi$ are equivalent in $\langle M,f\rangle$, with $f\in Reg$.

In subsection 2.4 we introduce the even finer class $Dec$ of {\em $\neg$-decreasing} regular associative  functions, that is,   $Dec\subset Reg^*$. A total ordering of $Sen(L)$ $<$ is $\neg$-decreasing if and only if  for all $\alpha$, $\beta$, $\alpha<\beta\Leftrightarrow \neg\beta<\neg\alpha$. $f$ is $\neg$-decreasing if and only if  $f=\min_<$ for some $\neg$-decreasing total ordering $<$.
The existence of $\neg$-decreasing  regular total orderings of $Sen(L)$ is shown and  a syntactic characterization of $\neg$-decreasingness is given.

In subsection 2.5 we consider the dual connective $\varphi\circ\psi:=\neg(\neg\varphi|\neg\psi)$ of $|$ and show that it commutes with $|$ if and only  the choice functions involved are $\neg$-decreasing.

Section 3 is devoted to the axiomatization of  Propositional Superposition Logic(s) (PLS). In the general section  we give axiomatizations for the logics based on the sets of choice functions ${\cal F}$,  $Reg$, $Reg^*$ and $Dec$. In general for every set $X\subseteq {\cal F}$ of choice functions and every set $K\subseteq Taut(X)$ of tautologies with respect to the truth notion $\models_X$, a logic ${\rm PLS}(X,K)$ is defined, whose axioms are those of PL plus $K$ and its semantics is the relation $\models_X$. Within ${\rm PLS}(X,K)$ $K$-consistency is defined and Soundness Theorem is proved for every logic ${\rm PLS}(X,K)$ with $K\subseteq Taut(X)$. Next  we introduce specific axiomatizations (by finitely many schemes of axioms) $K_0$, $K_1$, $K_2$, $K_3$ for the  truth relations   defined by the  classes ${\cal F}$, $Reg$, $Reg^*$ and $Dec$, respectively. The  logics ${\rm PLS}({\cal F},K_0)$, ${\rm PLS}(Reg,K_1)$, ${\rm PLS}(Reg^*,K_2)$,  ${\rm PLS}(Dec,K_3)$ are sound as a consequence of the previous general fact. There exists an essential difference between the axiomatization of ${\cal F}$, and those of the rest systems $Reg$, $Reg^*$ and $Dec$. The difference consists in  that $K_1$-$K_3$ contain an extra inference rule (besides Modus Ponens) because of which  the truth of the Deduction Theorem (DT) is open. This has serious effects on the completeness of the systems based on $K_1$-$K_3$. So we split the examination of completeness for ${\rm PLS}({\cal F},K_0)$ on the one hand and for  the rest systems on the other.

In subsection 3.1 we prove the (unconditional) completeness of the system ${\rm PLS}({\cal F},K_0)$.

In subsection  3.2 we examine completeness for the  logics  ${\rm PLS}(Reg,K_1)$, ${\rm PLS}(Reg^*,K_2)$ and ${\rm PLS}(Dec,K_3)$. The possible failure of DT makes it necessary to distinguish between two forms of completeness, CT1 and CT2, which in the lack of DT need not be equivalent. CT1 implies CT2 but the converse is open. Concerning the  systems $K_1$-$K_3$, we are seeking to prove  CT2 rather than CT1.  We show that these systems are {\em conditionally complete} in the sense that each of these systems is  CT2-complete if and only if each $K_i$ satisfies a certain extendibility property $cext(K_i)$ saying that every $K_i$-consistent set of sentences is extended to a $K_i$-consistent and complete set. This property is trivial for formal  systems $K$ satisfying  DT, but is  open for systems for which DT is open.  Assuming that $cext(K_i)$ is true, the proofs of CT2-completeness for the above logics are all variations of the proof of  completeness of ${\rm PLS}({\cal F},K_0)$. On the other hand failure of  $cext(K_i)$ implies the failure of CT2-completeness of the corresponding system.

In general, the proof of (CT2-)completeness of a logic ${\rm PLS}(X,K)$, with $K\subset Taut(X)$, goes roughly as follows: start with a $K$-consistent and complete set $\Sigma\subset Sen(L_s)$. To prove  it is $X$-verifiable, pick  $\Sigma_1=\Sigma\cap Sen(L)$. Then $\Sigma_1$ is a consistent and complete set in the sense of PL. So by completeness of the latter there exists a two-valued  assignment $M$ such that $M\models \Sigma_1$. Then in order to prove the $X$-verifiability of $\Sigma$, it suffices to  define a choice  function $f$ such that $f\in X$ and $\langle M,f\rangle \models\Sigma$.

Finally in section 4 we describe  briefly two  goals for future research, namely, (1) the goal to find alternative   semantics for the logics PLS, and (2) to develop a superposition extension of first-order logic (FOL) with an appropriate semantics and complete axiomatization.

\section{Semantics of  superposition propositional logic based on choice functions}

\subsection{Definitions and basic facts}

Let us fix a propositional language  $L=\{p_0,p_1,\ldots\}\cup\{\neg,\wedge\}$, where $p_i$ are  symbols of atomic sentences. The other connectives  $\vee$, $\rightarrow$, $\leftrightarrow$ are defined as usual in terms of the previous ones.\footnote{The functions $\overline{f}$ defined below are going to respect classical connectives, and hence classical equivalences, so it makes no difference if, e.g., we define $\varphi\rightarrow \psi$ as $\neg(\varphi\wedge \neg\psi)$ or $\neg(\neg\psi \wedge \phi)$.}  Let  $Sen(L)$ denote the set of sentences of $L$.  Throughout $M$ will denote  some  {\em truth assignment} for the sentences of $L$, that is, a mapping $M:Sen(L)\rightarrow \{0,1\}$ that is defined according to  the standard truth tables. For a given $\alpha\in Sen(L)$ we shall use the notation $M\models\alpha$, $M\models\neg\alpha$ instead  of $M(\alpha)=1$ and  $M(\alpha)=0$, respectively, for practical reasons. Namely, below we shall frequently refer to the truth of sentences  denoted $\overline{f}(\varphi)$,  so it would be more convenient to write $M\models\overline{f}(\varphi)$  than  $M(\overline{f}(\varphi))=1$.

Let $L_s=L\cup\{|\}$, where $|$ is a new primitive binary logical connective. The set of atomic sentences of $L_s$ are  identical to those of $L$, while the set of sentences of $L_s$, $Sen(L_s)$, is recursively defined along the  obvious steps: If $\varphi,\psi\in Sen(L_s)$, then $\varphi\wedge \psi$, $\varphi|\psi$, $\neg\varphi$ belong to $Sen(L_s)$.

\vskip 0.1in

{\bf Basic notational convention.} To keep track of whether we refer,  at each particular moment, to sentences of $L$ or  $L_s$, throughout the letters $\varphi$, $\psi$, $\sigma$ will denote general sentences  of  $L_s$, while the letters $\alpha$, $\beta$, $\gamma$ will denote sentences of  $L$ only. Also we often  refer to sentences of $L$  as ``classical''.

\vskip 0.1in

Throughout given a set $A$ we let
$$[A]^2=\{\{a,b\}:a,b\in A\}$$
denote the set of all 2-element and 1-element subsets of $A$.  We refer to the elements of $[A]^2$ as {\em pairs} of elements of $A$. A {\em choice function} for $[A]^2$ is as usual a mapping $f:[A]^2\rightarrow A$ such that $f(\{a,b\})\in \{a,b\}$ for every $\{a,b\}\in [A]^2$. To save brackets we write $f(a,b)$ instead of $f(\{a,b\})$. So in particular $f(a,b)=f(b,a)$ and $f(a,a)=a$.\footnote{The claim of the existence of a choice function for the set $[A]^2$, for every set $A$, is a weak form of the axiom of choice ($AC$), denoted ${\rm C}_2$ in \cite{Je73}. In general  for every $n\in\N$, ${\rm C}_n$ denotes the principle that every set of $n$-element sets has a choice function. The interested reader may consult  \cite[section 7.4]{Je73} for various interrelations between such principles, as well as with the Axiom of Choice for Finite Sets (saying that every nonempty set of nonempty  finite sets has a choice function). See also Remark \ref{R:2choice} below.}

\begin{Def} \label{D:choicefun}
{\em Given the language $L$, a choice function for $[Sen(L)]^2$, the set of pairs of sentences of $L$, will be referred to as a}  choice function for $L$.
\end{Def}
Let $${\cal F}(L)=\{f: \ \mbox{$f$ is a choice function for $L$}\}.$$
Throughout we shall write more simply ${\cal F}$ instead of ${\cal F}(L)$.
Below the letters $f,g$ will range  over elements of ${\cal F}$  unless otherwise stated. In particular for all $\alpha,\beta\in Sen(L)$, we write  $f(\alpha,\beta)$ instead of  $f(\{a,b\})$ so  $f(\alpha,\beta)=f(\beta,\alpha)$ and $f(\alpha,\alpha)=\alpha$.

\begin{Def} \label{D:collapse}
{\em Let $f$ be a choice function for $L$. Then $f$ generates a}  collapsing function {\em $\overline{f}:Sen(L_s)\rightarrow Sen(L)$ defined inductively as follows:

(i) $\overline{f}(\alpha)=\alpha$ for every $\alpha\in Sen(L)$.

(ii) $\overline{f}(\varphi\wedge \psi)=\overline{f}(\varphi)\wedge \overline{f}(\psi)$.

(iii) $\overline{f}(\neg \varphi)=\neg \ \overline{f}(\varphi)$.

(iv) $\overline{f}(\varphi|\psi)=f(\overline{f}(\varphi),\overline{f}(\psi))$.}
\end{Def}

\begin{Rems} \label{R:uniform}

{\em (i) Since the connectives $\vee$ and  $\rightarrow$ are  defined in terms of  $\neg$ and  $\wedge$, $\overline{f}$ commutes also with respect to them, that is,

\noindent $\overline{f}(\varphi\vee \psi)=\overline{f}(\varphi)\vee \overline{f}(\psi)$.

\noindent $\overline{f}(\varphi\rightarrow \psi)=\overline{f}(\varphi)\rightarrow \overline{f}(\psi)$.

(ii) The crucial clause of the definition is of course (iv). It says that for any sentences $\varphi$, $\psi$, $\overline{f}(\varphi|\psi)$ is a} choice {\em from the set $\{\overline{f}(\varphi),\overline{f}(\psi)\}$. In particular, for  classical sentences $\alpha$, $\beta$ we have
\begin{equation} \label{E:crucial}
\overline{f}(\alpha|\beta)=f(\alpha,\beta),
\end{equation}
}
\end{Rems}

\begin{Def} \label{D:s-truth}
{\em (Main Truth definition). Let $M$ be a truth assignment for $L$, $f$ a choice function for $L$ and $\overline{f}:Sen(L_s)\rightarrow Sen(L)$ be  the corresponding collapsing function. The} truth relation $\models_s$ between the pair $\langle M,f\rangle$ and a sentence $\varphi$ of $L_s$ {\em is defined as follows:
$$\langle M,f\rangle\models_s \varphi \ \mbox{iff} \ M\models \  \overline{f}(\varphi).$$
More generally, for a set $\Sigma\subset Sen(L_s)$  we write $\langle M,f\rangle\models_s\Sigma$, if $\langle M,f\rangle\models_s \varphi$ for every $\varphi\in\Sigma$.
}
\end{Def}

The following facts are easy consequences of the preceding definitions.

\begin{Fac} \label{F:extend}
(i) The truth relation $\models_s$ extends the Boolean  one $\models$, that is for every $\alpha\in Sen(L)$, and every $\langle M,f\rangle$, $\langle M,f\rangle\models_s\alpha \Leftrightarrow \ M\models\alpha$.

(ii) $\models_s$ is a  bivalent notion of truth, that is for every $\langle M,f\rangle$ and every sentence $\varphi$,  either $\langle M,f\rangle\models_s \varphi$ or $\langle M,f\rangle\models_s \neg\varphi$.

(iii) For every sentence  $\varphi$ of $L_s$, every structure $M$ and  every collapsing function  $\overline{f}$, $\langle M,f\rangle\models_s\varphi|\varphi$ if and only if  $\langle M,f\rangle\models_s\varphi$.

(iv)  For all $\varphi,\psi\in Sen(L_s)$,  $M$ and $f$, $\langle M,f\rangle\models_s\varphi|\psi$ if and only if  $\langle M,f\rangle\models_s\psi|\varphi$.
\end{Fac}

{\em Proof.} (i) Immediate from the fact that by clause (i) of \ref{D:collapse}, $\overline{f}(\alpha)=\alpha$ for every sentence $\alpha\in Sen(L)$. Thus $\langle M,f\rangle\models_s\alpha$ if and only if  $M\models\alpha$.

(ii) Let $\langle M,f\rangle\not\models_s\varphi$. Then $M\not\models \overline{f}(\varphi)$, that is, $M\models \neg\overline{f}(\varphi)$. By clause  (iii) of \ref{D:collapse}, $\neg\overline{f}(\varphi)=\overline{f}(\neg\varphi)$, so $\langle M,f\rangle\not\models_s\varphi$ implies $M\models \overline{f}(\neg \varphi)$, or $\langle M,f\rangle\models_s\neg\varphi$.

(iii): By definition \ref{D:choicefun}, $f(\alpha,\alpha)=\alpha$, for every $\alpha$. Therefore
$\langle M,f\rangle\models_s\varphi|\varphi \Leftrightarrow {\cal M}\models f(\overline{f}(\varphi),\overline{f}(\varphi)) \Leftrightarrow {\cal M}\models\overline{f}(\varphi)\Leftrightarrow \langle M,f\rangle\models_s\varphi$.

(iv): By  \ref{D:choicefun} again $f(\alpha,\beta)=f(\beta,\alpha)$. So
$\langle M,f\rangle\models_s\varphi|\psi \Leftrightarrow M\models f(\overline{f}(\varphi),\overline{f}(\psi)) \Leftrightarrow M\models f(\overline{f}(\psi),\overline{f}(\varphi))\Leftrightarrow \langle M,f\rangle\models_s\psi|\varphi$. \telos

\vskip 0.2in

Let $\Sigma\models\alpha$, $\alpha\models\beta$,  (for $\Sigma\subset Sen(L$)), and $\alpha\sim\beta$ denote the classical logical consequence and logical equivalence relations, respectively, for classical sentences. These are extended to the relations $\varphi\models_s\psi$, $\Sigma\models_s\varphi$ (for $\Sigma\subset Sen(L_s$)),  and $\varphi\sim_s\psi$ for $L_s$-sentences as follows.

\begin{Def} \label{D:logiceq}
{\em Let $\Sigma\subset Sen(L)$, $\varphi,\psi\in Sen(L_s)$. We say that $\varphi$ is an} s-logical consequence {\em of $\Sigma$, denoted $\Sigma\models_s\varphi$, if for every  structure $\langle M,f\rangle$, $\langle M,f\rangle\models_s\Sigma$ implies $\langle M,f\rangle\models_s\varphi$. In particular we write $\varphi\models_s\psi$ instead of $\{\varphi\}\models_s\psi$.
We say that $\varphi$ and $\psi$ are} s-logically  equivalent, {\em  denoted $\varphi\sim_s\psi$, if for every  $\langle M,f\rangle$,
$$\langle M,f\rangle\models_s\varphi \ \Leftrightarrow \ \langle M,f\rangle\models_s\psi.$$
Finally, $\varphi$ is an} s-tautology, {\em denoted $\models_s\varphi$, if $\langle M,f\rangle\models_s \varphi$ for every $\langle M,f\rangle$.
}
\end{Def}

\begin{Fac} \label{F:tautology}
(i) $\varphi\models_s\psi$ if and only if  $\models_s\varphi\rightarrow \psi$.

(ii) $\varphi\sim_s\psi$ if and only if  $\models_s\varphi\leftrightarrow \psi$.

(iii) For $\alpha,\beta\in Sen(L)$, $\alpha\models_s\beta$ if and only if  $\alpha\models\beta$ and $\alpha\sim_s\beta$ if and only if  $\alpha\sim\beta$.

(iv) $\varphi\sim_s\psi$ if and only if  for all choice functions $f$, $\overline{f}(\varphi)\sim \overline{f}(\psi)$.

(v) Let $\alpha(p_1,\ldots,p_n)$ be a  sentence of $L$, made up by the atomic sentences $p_1,\ldots,p_n$, let $\psi_1,\ldots,\psi_n$ be any sentences of $L_s$ and let $\alpha(\psi_1,\ldots,\psi_n)$ be the sentence resulting from $\alpha$ if we replace each $p_i$ by $\psi_i$. Then: $$\models\alpha(p_1,\ldots,p_n) \ \Rightarrow \ \models_s\alpha(\psi_1,\ldots,\psi_n).$$

(vi) For all $\varphi,\psi$, $\varphi|\varphi\sim_s\varphi$ and $\varphi|\psi\sim_s\psi|\varphi$.

(vii) If $\varphi\sim_s\psi$, then $\varphi|\psi\sim_s\varphi$.

\end{Fac}

{\em Proof.} (i): Let $\varphi\models_s\psi$. It means that for every $\langle M,f\rangle$, $\langle M,f\rangle\models_s\varphi$ implies $\langle M,f\rangle\models_s\psi$. Equivalently,  $M\models\overline{f}(\varphi)$ implies $M\models\overline{f}(\psi)$, or $M\models(\overline{f}(\varphi)\rightarrow\overline{f}(\psi))$, or $M\models\overline{f}(\varphi\rightarrow\psi)$. It means that for every $M$ and every $\overline{f}$, $\langle M,f\rangle\models_s\varphi\rightarrow\psi$. Thus $\models_s\varphi\rightarrow\psi$. The converse is similar.

(ii) and  (iii) follow from (i).

(iv): Note that $\varphi\sim_s\psi$ holds if and only if  for all $\langle M,f\rangle$, $\langle M,f\rangle\models_s\varphi$ if and only if  $\langle M,f\rangle\models_s\psi$, or equivalently, $M\models\overline{f}(\varphi)$ if and only if  $M\models\overline{f}(\psi)$. But this means that for every $f$, $\overline{f}(\varphi)\sim \overline{f}(\psi)$.

(v): Suppose $\models\alpha(p_1,\ldots,p_n)$. For any choice function   $f$, clearly $$\overline{f}(\alpha(\psi_1,\ldots,\psi_n))=
\alpha(\overline{f}(\psi_1),\ldots,\overline{f}(\psi_n)),$$
since $\alpha$ is classical and $\overline{f}$ commutes with standard connectives. Moreover $\models \alpha(\overline{f}(\psi_1),\ldots,\overline{f}\psi_n)$, since by assumption $\models\alpha(p_1,\ldots,p_n)$ and $\overline{f}(\psi_i)$ are standard sentences. Thus $M\models\alpha(\overline{f}(\psi_1),\ldots,\overline{f}(\psi_n))$, for every $M$, or $M\models\overline{f}(\alpha(\psi_1,\ldots,\psi_n))$. It means that
$\langle M,f\rangle\models_s\alpha(\psi_1,\ldots,\psi_n)$ for every structure $\langle M,f\rangle$, so  $\models_s\alpha(\psi_1,\ldots,\psi_n)$.

(vi) This follows immediately from clauses (iii) and (iv) of Fact \ref{F:extend}.

(vii) Let $\varphi\sim_s\psi$ and let  $\langle M,f\rangle\models \varphi|\psi$. Then  $M\models f(\overline{f}(\varphi),\overline{f}(\psi))$. By clause (iv) above, $\overline{f}(\varphi)\sim\overline{f}(\psi)$ since $\varphi\sim_s\psi$. Therefore whatever the choice of $f$ would be between $\overline{f}(\varphi)$ and $\overline{f}(\psi)$, we shall have $M\models\overline{f}(\varphi)$. Thus $\langle M,f\rangle\models_s\varphi$.
\telos

\vskip 0.2in

The following {\em interpolation  property} of the new semantics is perhaps  the most striking one. Notice  that it holds for the {\em general  choice functions,}  not requiring any of the additional conditions to be  considered in the subsequent sections.

\begin{Thm} \label{T:interpol}
For all $\varphi,\psi\in Sen(L_s)$,
$$\varphi\wedge \psi \models_s  \varphi|\psi \models_s  \varphi\vee \psi,$$
while in general
$$\varphi\vee \psi \not\models_s\varphi|\psi \not\models_s \varphi\wedge \psi.$$
\end{Thm}

{\em Proof.}  Assume $\langle M,f\rangle\models_s \varphi\wedge\psi$. Then $M\models \overline{f}(\varphi)\wedge \overline{f}(\psi)$, that is, $M\models \overline{f}(\varphi)$ and  $M\models \overline{f}(\psi)$. But then, whatever $f$ would choose from  $\{\overline{f}(\varphi),\overline{f}(\psi)\}$, it  would be true in $M$, that is, $M\models f(\overline{f}(\varphi),\overline{f}(\psi))$. This exactly  means that $\langle M,f\rangle\models_s \varphi|\psi$. Therefore $\varphi\wedge \psi \models_s  \varphi|\psi$.

On the other hand, if  $\langle M,f\rangle\models_s \varphi|\psi$ then  $M\models f(\overline{f}(\varphi),\overline{f}(\psi))$. If  $f(\overline{f}(\varphi),\overline{f}(\psi))=\overline{f}(\varphi)$, then $M\models \overline{f}(\varphi)$.  If $f(\overline{f}(\varphi),\overline{f}(\psi))=\overline{f}(\psi)$, then $M\models \overline{f}(\psi)$. So either $M\models \overline{f}(\varphi)$ or $M\models \overline{f}(\psi)$. Therefore  $M\models \overline{f}(\varphi)\vee\overline{f}(\psi)$. But clearly $\overline{f}(\varphi)\vee\overline{f}(\psi)=\overline{f}(\varphi\vee \psi)$, since $\overline{f}$ commutes with all  standard connectives. Thus $M\models \overline{f}(\varphi\vee \psi)$, or equivalently, $\langle M,f\rangle\models_s \varphi\vee\psi$. Therefore $\varphi|\psi \models_s  \varphi\vee\psi$.

To see that the converse relations are false, pick $\alpha\in Sen(L)$ and a truth assignment $M$ such that  $M\models\alpha$. Pick also a choice function for $L$ such that  $f(\alpha,\neg\alpha)=\neg\alpha$. Since  $\overline{f}(\alpha\vee\neg\alpha)=\alpha\vee\neg\alpha$,  $\langle M,f\rangle\models_s\alpha\vee \neg\alpha$. On the other hand, $M\not\models\neg\alpha$ implies $M\not\models f(\alpha,\neg\alpha)$, thus   $\langle M,f\rangle\not\models_s \alpha|\neg\alpha$. Therefore  $\alpha\vee \neg\alpha \not\models_s\alpha|\neg\alpha$. Similarly, if $M$, $\alpha$ are as before, but we take a choice function $g$ such that $g(\alpha,\neg\alpha)=\alpha$, then $\langle M,g\rangle\models_s \alpha|\neg\alpha$, while  $\langle M,g\rangle\not\models_s \alpha\wedge\neg\alpha$. So $\alpha|\neg\alpha \not\models_s \alpha\wedge \neg\alpha$. \telos

\vskip 0.2in

\begin{Cor} \label{C:Sch}
If $\alpha$ is neither a  tautology nor a contradiction, then  $\alpha|\neg\alpha$ is neither an $s$-tautology nor an $s$-contradiction.
\end{Cor}

{\em Proof.} If  $\alpha$ is  as stated, then by the proof of Theorem \ref{T:interpol} $\alpha|\neg\alpha$ is strictly interpolated between $\alpha\wedge\neg\alpha$ and $\alpha\vee\neg\alpha$.  \telos

\vskip 0.2in

In the semantics $\models_s$ used above, arbitrary  choice functions for $L$ are allowed to participate. This practically  means that for any  pair  $\{\alpha,\beta\}$, $f$ may pick an element from $\{\alpha,\beta\}$ quite randomly, e.g. by tossing a coin. However, if we want  $\models_s$ to support additional properties of $|$,  we must refine $\models_s$ by imposing  extra  conditions to the choice functions. Such a refinement can be defined in a general manner as follows.

\begin{Def} \label{D:X-notion}
{\em For every $\emptyset\neq X\subseteq {\cal F}$, define the} $X$-logical consequence relation $\models_X$ {\em and the} $X$-logical equivalence relation $\sim_X$ {\em  as follows:  $\Sigma\models_X\varphi$  if and only if  for every truth assignment $M$ for $L$ and every $f\in X$, $$\langle M,f\rangle\models_s\Sigma \ \Rightarrow \langle M,f\rangle\models_s\varphi.$$
Also  $\varphi\sim_X\psi$ if and only if  $\varphi\models_X\psi$ and $\psi\models_X\varphi$.}
\end{Def}
[The purpose of condition $X\neq \emptyset$  is to  block trivialities. For if $X=\emptyset$, we vacuously have  $\varphi\models_\emptyset\psi$ and $\varphi\sim_\emptyset\psi$ for all $\varphi,\psi\in Sen(L_s)$. So  all sets $X$, $Y\subseteq {\cal F}$ referred to below are  assumed to be $\neq\emptyset$.]

Using the  above  notation, the relations $\models_s$ and $\sim_s$ are  alternatively written  $\models_{\cal F}$, $\sim_{\cal F}$, respectively.

The following simple fact reduces $\sim_X$ to the standard $\sim$.

\begin{Lem} \label{L:simplereduction}
For every $X\subseteq {\cal F}$, and any $\varphi,\psi\in Sen(L_s)$,
$$\varphi\sim_X\psi \ \Leftrightarrow (\forall f\in X)(\overline{f}(\varphi)\sim \overline{f}(\psi)).$$
\end{Lem}

{\em Proof.} By  definition, $\varphi\sim_X\psi$ if for every $M$ and every $f\in X$, $$\langle M,f\rangle\models_s\varphi \ \Leftrightarrow \ \langle M,f\rangle\models_s\psi,$$ or, equivalently, if for all $M$ and  $f\in X$, $$M\models\overline{f}(\varphi) \ \Leftrightarrow \ M\models\overline{f}(\psi).$$
The latter is true if and only if   for all $f\in X$, $\overline{f}(\varphi)\sim \overline{f}(\psi)$. \telos

\vskip 0.2in

The next properties are  easy to verify.

\begin{Fac} \label{F:tautolol}
For every $X,Y\subseteq {\cal F}$:

(i) $\varphi\models_X\psi$ if and only if  $\models_X\varphi\rightarrow \psi$.

(ii) $\varphi\sim_X\psi$ if and only if  $\models_X\varphi\leftrightarrow \psi$.

(iii) If $X\subseteq Y$, then $\models_Y\subseteq \models_X$ and $\sim_Y\subseteq \sim_X$.

(iv) The restriction of $\sim_X$ to classical sentences coincides with $\sim$, that is, for all $\alpha,\beta\in Sen(L)$, $$\alpha\sim_X\beta \ \Leftrightarrow \ \alpha\sim\beta.$$
\end{Fac}

Before closing this section  I should  give credit to \cite{Hum11}  for some notions introduced above.
It was not until one of the referees drew my attention to \cite{Hum11}, when I learned (with surprise) that the notion of  choice function for pairs of formulas, and, essentially, the germ of the satisfaction  relation defined in \ref{D:s-truth} above, were not entirely new but had already been defined independently with some striking similarities in the style of presentation. In fact in Example 3.24.14, p. 479, of \cite{Hum11} we read:
\begin{quote}
``By a {\em  pair selection function }  on a set $U$ we mean a function $f$ such that for all $a,b\in U$, $f(\{a,b\})\in \{a,b\}$. We write $f(a,b)$ for `$f(\{a,b\}$' and include the possibility that $a=b$ in which case $f(a,b)=a=b$. (...) A pair selection function is accordingly a commutative idempotent binary operation which is in addition a {\em quasi-projection} or a {\em conservative} operation, meaning that its value for a given pair of arguments is always one of those arguments. For the current application consider $f$ as a pair selection function on the set of formulas of the language generated from the actual stock of propositional variables with the aid of the binary connective $\circ$. Consider the gcr (=generalized consequence relation) determined by the class of all valuations $v$ satisfying the condition that for some pair selection function $f$ we have: For all formulas $A,B$, $v(A\circ B)=v(f(A,B))$. Then, if $\succ$ denotes this gcr, it satisfies the rules: (I) $A,B\succ A\circ B$, (II) $A\circ B\succ A,B$ and (IV) $A\circ B \succ B\circ A$.''
\end{quote}

Note that rules (I) and (II) are essentially the ``interpolation property''  of Theorem \ref{T:interpol}, while rule (IV) is the commutativity property (vi) of Fact \ref{F:tautology}.

In the next section we consider a first natural subclass $X\subset{\cal F}$, the class of {\em associative} choice functions. These are precisely the functions that support the associativity property of the connective $|$. Clearly associativity is a highly desirable property  from an {\em algebraic} point of view. However, as one of the referees interestingly observed at this point,  we must distinguish between what is  algebraically desirable  and what is  quantum mechanically desirable, i.e., close to the real  behavior of a quantum system. In his view, classes of choice functions with not very  attractive and smooth properties might also deserve to  be isolated and scrutinized.

\subsection{Associative choice functions}

By clause (vi) of Fact \ref{F:tautology},  $\varphi|\varphi\sim_s\varphi$ and  $\varphi|\psi\sim_s\psi|\varphi$. These two properties, idempotence  and commutativity up to logical equivalence,  are in accord with the intended intuitive meaning of the operation $|$.  Another  desirable property that is in accord with the meaning of  $|$ is  {\em associativity},  that is, the logical equivalence of  $(\varphi|\psi)|\sigma$ and $\varphi|(\psi|\sigma)$. Is it true with respect to $\sim_s$? The answer is: not in general.
In order to ensure it we  need to impose a certain condition on the choice functions.   The specific condition  does not depend on the nature of  elements of $Sen(L)$, so we prove it below in a general setting.

\vskip 0.2in

Let $A$ be an infinite set and let  $f:[A]^2\rightarrow A$ be a choice for  pairs of elements of $A$. One might extend it to the set $[A]^3$, of nonempty sets  with at most $3$ elements, by setting, say,
$$f(a,b,c):=f(\{a,b,c\})=f(f(a,b),c).$$
But this does not guarantee that $f(a,b,c)=f(b,c,a)=f(c,b,a)$, etc, as it would be obviously required, unless
$f(f(a,b),c)=f(a,f(b,c))$ for all $a,b,c$. This is exactly the required condition.

\begin{Def} \label{D:associative}
{\em Let $f$ be a choice function for $[A]^2$. $f$ is said to be} associative {\em if for all $a,b,c\in A$, $$f(f(a,b),c)=f(a, f(b,c)).\footnote{
If we write $a\star b$ instead of $f(a,b)$, then  the condition  $f(f(a,b),c)=f(a, f(b,c))$ is rewritten $(a\star b)\star c=a\star(b\star c)$, which justifies the term ``associative''.}$$}
\end{Def}

We show below that associative choice functions on $[A]^2$ are, essentially,   the functions $\min_{<}$, where $<$ is a total ordering of $A$. I do not know if the next theorem is new or a known result. In any case I couldn't find a proof in the current bibliography.

\begin{Thm} \label{T:LO}
(i) If $<$ is a total ordering on $A$, then the mapping  $\min_{<}(a,b)$ from $[A]^2$ into $A$ is associative.

(ii) Conversely, if $f:[A]^2\rightarrow A$ is an associative choice function, then it defines  a total ordering $<$ on $A$ such that for all $a,b\in A$, $f(a,b)=\min_<(a,b)$.
\end{Thm}

{\em Proof.} (i) Let $<$ be a total ordering of $A$. Let $Fin(A)$ denote the set of all nonempty finite subsets of $A$ and let $\min_<$ be the function picking the  $<$-least element of $x$ for every $x\in Fin(A)$. Let us write $\min$ instead of  $\min_<$. Obviously $\min$  is a   choice function for $Fin(A)$. In particular, for all $a,b,c\in A$,
$$\min(a,b,c)=\min(\min(a,b),c)=\min(a,\min(b,c)).$$
Thus $\min$ restricted to $[A]^2$ is associative.

(ii) Let $f:[A]^2\rightarrow A$ be an associative choice function.
Define the relation $<$ on $A$ as follows:
For any $a,b\in A$, let   $a<b$ if and only if  $a\neq b$ and $f(a,b)=a$. Obviously $<$ is  total and anti-reflexive  (that is, $a\not<a$ for every $a\in A$). Thus in order for $<$ to be a total ordering it suffices to be also  transitive. Let $a<b$ and $b<c$. We show that $a<c$. By the assumptions, we have  $a\neq b$,  $f(a,b)=a$, $b\neq c$ and $f(b,c)=b$. It follows from them that  $a\neq c$, for otherwise $b=f(b,c)=f(b,a)=f(a,b)=a$, a contradiction.
It remains to show that $f(a,c)=a$. By associativity and commutativity of $f$,   $f(a,f(b,c))=f(b,f(a,c))$. Since  $f(a,f(b,c))=f(a,b)=a$, it follows that
$f(b,f(a,c))=a$ too. If $f(a,c)=c$, then we would have
$f(b,f(a,c))=f(b,c)=b\neq a$,  a contradiction. Therefore  $f(a,c)=a$ and we are done.  Thus $<$ is a total ordering of $A$, and by definition  $f(a,b)=\min_{<}(a,b)$, for all $a,b\in A$. \telos

\vskip 0.2in

As an immediate corollary of Theorem \ref{T:LO} we obtain the following.

\begin{Cor} \label{C:2extends}
If $f:[A]^2\rightarrow A$ is an associative choice function, then it defines uniquely a total ordering $<$ of $A$ such that $f(a,b)=\min_<(a,b)$. Therefore $f$ extends uniquely to the choice function  $f^+:Fin(A)\rightarrow A$, such that for every $x\in Fin(A)$, $f^+(x)=\min_<(x)$.  Thus $f=f^+\restr [A]^2$.
\end{Cor}

In view of the preceding Corollary, we can without serious loss of precision  identify an  associative choice function $f$ on  $[A]^2$ with the generated choice function $f^+$ on the entire $Fin(A)$,  and write  $f=\min_<$ instead of $f^+=\min_<$, where $<$ is the ordering defined by $f$.

\begin{Rem} \label{R:2choice}
{\em From a set-theoretical point of view, the existence of an associative function is a much stronger statement than the existence of a simple choice function for $[A]^2$. As noticed in footnote 9 the latter is identical to the principle ${\rm C}_2$. On the other hand, it follows from   \ref{T:LO} and \ref{C:2extends} that the existence of an associative choice function for $[A]^2$, for  every set $A$, is equivalent to the existence of a total ordering on $A$, i.e., to the} Ordering Principle {\em saying that ``Every set can be totally ordered'', which is much stronger than ${\rm C}_2$. Specifically, it was  shown in \cite{La64} that the Ordering Principle is strictly stronger than the Axiom of Choice for Finite Sets. The latter is in turn strictly stronger than the conjunction of all axioms ${\rm C}_n$, for  $n\geq 2$ (see \cite[Theorem 7.11]{Je73}).}
\end{Rem}

Let us now return to the set $Sen(L)$ of sentences of $L$. In particular a  choice function $f$ for $L$ is said to be {\em associative} if  for all $\alpha,\beta,\gamma\in Sen(L)$,
\begin{equation} \label{E:assphi}
f(f(\alpha,\beta),\gamma)=f(\alpha, f(\beta,\gamma)).
\end{equation}
We often call also the pair $\langle M,f\rangle$ {\em associative} if $f$ is associative. As an immediate consequence of Theorem \ref{T:LO}, Corollary \ref{C:2extends} and the comments following the latter we have  the following.

\begin{Cor} \label{C:addasso}
A choice function $f$ for $L$ is associative if and only if  there is a total ordering $<$ of $Sen(L)$ such that $f=\min_<$.
\end{Cor}

Let
$$Asso=\{f\in {\cal F}: f \ \mbox{is  associative }\},$$
denote the set of all associative choice functions for $L$. We shall also  abbreviate
the logical consequence relation $\models_{Asso}$ and the logical equivalence relation  $\sim_{Asso}$ that are induced by  $Asso$ (see definition \ref{D:X-notion} of the previous section), by $\models_{Asso}$ and $\sim_{Asso}$, respectively. From the general facts \ref{F:tautolol} (iii) we immediately obtain the following.

\begin{Fac} \label{F:stronass}
For all $\varphi,\psi,\Sigma$,

(i) $\Sigma\models_s\varphi \ \Rightarrow \ \Sigma\models_{Asso}\varphi$.

(ii) $\varphi\sim_s\psi \ \Rightarrow \ \varphi\sim_{Asso}\psi$.
\end{Fac}

\noindent It is easy to verify that the arrows in the preceding Fact cannot in general be  reversible. The main consequence of associativity  is the following.

\begin{Thm} \label{T:assoc}
Let $X\subseteq {\cal F}$ be a class of choice functions. If $X\subseteq Asso$, then  $|$ is associative with respect to the truth notion $\models_X$,  that is, for all $\varphi$, $\psi$, $\sigma$, $\varphi|(\psi|\sigma)\sim_X (\varphi|\psi)|\sigma$.
\end{Thm}

{\em Proof.}  Let $X\subseteq Asso$. It suffices to show that for every $M$ and every $f\in X$, and any sentences  $\varphi$, $\psi$, $\sigma$ of $L_s$,  $\langle M,f\rangle$, $$\langle M,f\rangle\models_s(\varphi|\psi)|\sigma \ \mbox{iff} \  \langle M,f\rangle\models_s\varphi|(\psi|\sigma).$$
Fix some $M$ and some $f\in X$.
By  definition  we have:
$$\langle M,f\rangle\models_s(\varphi|\psi)|\sigma \Leftrightarrow M\models f(\overline{f}(\varphi|\psi),\overline{f}(\sigma)) \Leftrightarrow M\models f(f(\overline{f}(\varphi), \overline{f}(\psi)),\overline{f}(\sigma)).$$
By assumption $f\in Asso$, so by the associativity property (\ref{E:assphi})
$$f(f(\overline{f}(\varphi), \overline{f}(\psi)),\overline{f}(\sigma))=
f(\overline{f}(\varphi),f(\overline{f}(\psi),\overline{f}(\sigma))).$$
Thus,
$$\langle M,f\rangle\models_s(\varphi|\psi)|\sigma
 \Leftrightarrow M\models f(\overline{f}(\varphi),\overline{f}(\psi|\sigma)) \Leftrightarrow \langle M,f\rangle\models_s\varphi|(\psi|\sigma).$$ \telos

\vskip 0.2in

If we slightly weaken the property of associativity, the converse of \ref{T:assoc} holds too.

\begin{Def} \label{D:essass}
{\em Let us call a choice function for $L$} essentially associative, {\em if  (\ref{E:assphi}) holds with $\sim$ in place  of $=$, that is, for all $\alpha,\beta,\gamma\in Sen(L)$,
\begin{equation} \label{E:essassphi}
f(f(\alpha,\beta),\gamma)\sim f(\alpha, f(\beta,\gamma)).
\end{equation}
}
\end{Def}

Let $Asso'$ denote the class of essentially associative choice functions.

\begin{Thm} \label{T:asso'}
If $X\subseteq {\cal F}$ and $|$ is associative with respect to $\models_X$, then $X\subseteq Asso'$.
\end{Thm}

{\em Proof.} Let  $X\not\subseteq Asso'$. We have to show that $|$ is not associative with respect to $\models_X$. Pick $f\in X-Asso'$. It suffices to find $M$ and $\alpha$, $\beta$, $\gamma$ in $Sen(L)$ such that
$$\langle M,f\rangle\models_s(\alpha|\beta)|\gamma\not\Leftrightarrow \langle M,f\rangle\models_s\alpha|(\beta|\gamma).$$
Since $f$ is not essentially associative, there are $\alpha$, $\beta$, $\gamma$ in $Sen(L)$, such that $f(f(\alpha,\beta),\gamma)\not\sim f(\alpha,f(\beta,\gamma))$. Without loss of generality we may assume that
$$f(f(\alpha,\beta),\gamma)=\gamma\not\sim \alpha=f(\alpha,f(\beta,\gamma)).$$
Since $\alpha\not\sim\gamma$ there is  $M$ such that $M\models\alpha\wedge\neg\gamma$ or $M\models\neg\alpha\wedge\gamma$. Without loss of generality assume that the first is the case.   Then $M\models\alpha=f(\alpha,f(\beta,\gamma))$, so  $M\models \overline{f}(\alpha|(\beta|\gamma))$, therefore $\langle M,f\rangle\models_s\alpha|(\beta|\gamma)$. On the other hand $M\not\models\gamma=f(f(\alpha,\beta),\gamma)$, which implies $\langle M,f\rangle\not\models_s(\alpha|\beta)|\gamma$. This proves the theorem. \telos

\vskip 0.1in

Obviously $Asso\subseteq Asso'$. Are the two classes distinct? The answer is yes, but the functions in $Asso'-Asso$  behave {\em non-associatively} only on sentences $\alpha,\beta,\gamma$ such that $\alpha\sim\beta\sim\gamma$. To be precise, let us say that a  triple of sentences $\alpha,\beta,\gamma$ {\em witnesses non-associativity} of $f$, if  $f(f(\alpha,\beta),\gamma)\neq f(\alpha,f(\beta,\gamma))$, or $f(f(\alpha,\beta),\gamma)\neq f(\beta,f(\alpha,\gamma))$, or $f(f(\alpha,\gamma),\beta)\neq f(\alpha,f(\beta,\gamma))$.  Then the following holds.

\begin{Lem} \label{L:dist}
(i) $Asso\varsubsetneq Asso'$.

(ii) If  $f\in Asso'$ and $\alpha,\beta,\gamma$ are sentences such that $\alpha\not\sim \beta$, $\beta\not\sim\gamma$ and $\alpha\not\sim\gamma$, then $f$ is associative on $\alpha,\beta,\gamma$, i.e.,
$f(f(\alpha,\beta),\gamma)= f(\alpha,f(\beta,\gamma))=f(f(\alpha,\gamma),\beta)$.

(iii) If $f\in Asso'-Asso$, and $\alpha,\beta,\gamma$ witness the non-associativity of $f$, then  $\alpha,\beta,\gamma$ are all distinct, and  besides  $f(\alpha,\beta)$, $f(\alpha,\gamma)$, $f(\beta,\gamma)$ are all distinct.

(iv) Therefore, if $f\in Asso'-Asso$ and $\alpha,\beta,\gamma$ witness the non-associativity of $f$, then $\alpha\sim\beta\sim\gamma$.

(v) Further, if $f\in Asso'-Asso$, then $f$ is associative on every triple $\alpha$, $\beta$, $\gamma$ such that $\alpha\sim\beta\not\sim\gamma$.
\end{Lem}

{\em Proof.} (i) Let $\alpha\sim\beta\sim\gamma$, while all $\alpha,\beta,\gamma$ are distinct. Let  $f\in {\cal F}$ be such that $f(\alpha,\beta)=\beta$, $f(\alpha,\gamma)=\alpha$, $f(\beta,\gamma)=\gamma$. Then obviously all $f(f(\alpha,\beta),\gamma)$, $f(f(\alpha,\gamma),\beta)$, $f(f(\beta,\gamma),\alpha)$ are equivalent, so $f\in Asso'$. On the other hand, for example, $f(f(\alpha,\beta),\gamma)\neq f(\alpha,f(\beta,\gamma))$, so $f\notin Asso$.

(ii) If $\alpha\not\sim \beta$, $\beta\not\sim\gamma$, $\alpha\not\sim\gamma$, then it cannot be $f(f(\alpha,\beta),\gamma)\sim f(\alpha,f(\beta,\gamma))$ unless $f(f(\alpha,\beta),\gamma)=f(\alpha,f(\beta,\gamma))$.

(iii) It is easy to see that for every $f\in {\cal F}$ and every $\alpha,\beta$, $f(f(\alpha,\beta),\alpha)=f(f(\alpha,\beta),\beta)$. This shows that if any two elements of a triple $\alpha,\beta,\gamma$ are equal, this triple cannot witness the non-associativity of any function. Let $f\in Asso'-Asso$ and suppose  $\alpha,\beta,\gamma$ witness the non-associativity of $f$. We have just seen that they are all distinct. We show that $f(\alpha,\beta)$, $f(\alpha,\gamma)$, $f(\beta,\gamma)$ are  distinct too. Indeed assume that  two of the values
$f(\alpha,\beta)$, $f(\alpha,\gamma)$, $f(\beta,\gamma)$ are identical. It will follow that
\begin{equation} \label{E:wit}
f(f(\alpha,\beta),\gamma)= f(\alpha,f(\beta,\gamma))=f(f(\alpha,\gamma),\beta),
\end{equation}
which contradicts  the fact that
$\alpha,\beta,\gamma$    witness the non-associativity of $f$. Assume without loss of generality that $f(\alpha,\beta)=f(\beta,\gamma)$. Since $f$ is a choice function and $\alpha$, $\beta$, $\gamma$ are distinct,  necessarily $f(\alpha,\beta)=f(\beta,\gamma)=\beta$.  Therefore
$$f(f(\alpha,\beta),\gamma)=f(\beta,\gamma)=f(\alpha,\beta)
=f(\alpha,f(\beta,\gamma))=\beta.$$
So two members of (\ref{E:wit}) are equal. As to the third one, observe that  $f(\alpha,\gamma)$ is either $\alpha$ or $\gamma$. In both cases  $f(f(\alpha,\gamma),\beta)=\beta$, as required.

(iv) Let $f\in Asso'-Asso$ and let $\alpha,\beta,\gamma$ witness the non-associativity of $f$. By (iii) above, $f(\alpha,\beta)$, $f(\alpha,\gamma)$, $f(\beta,\gamma)$ take up all the values $\alpha,\beta,\gamma$, and therefore so do $f(f(\alpha,\beta),\gamma)$,  $f(\alpha,f(\beta,\gamma))$, $f(f(\alpha,\gamma),\beta)$. But since $f\in Asso'$, the latter are all logically equivalent. Therefore $\alpha\sim\beta\sim\gamma$.

(v) If  $f\in Asso'-Asso$, $\alpha\sim\beta\not\sim\gamma$ and $f$ were not  associative on $\alpha,\beta,\gamma$, the latter triple would witness the non-associativity of $f$, so, by (iv),  $\alpha\sim\beta\sim\gamma$. A contradiction. \telos

\vskip 0.2in

It follows from the preceding Lemma that every $f\in Asso'$ defines essentially an associative choice function (and hence a total ordering) for  the set of pairs of elements of $Sen(L)/\!\!\sim=\{[\alpha]:\alpha\in Sen\}$  rather than $Sen(L)$.

By  Facts \ref{F:extend}, \ref{F:stronass} and Theorem  \ref{T:assoc}, we obtain  the following.

\begin{Cor} \label{C:allprop}
The operation $|$ is idempotent, commutative and associative with respect to $\sim_{Asso}$. That is:

(i) $\varphi|\varphi\sim_{Asso} \varphi$.

(ii)  $\varphi|\psi\sim_{Asso} \psi|\varphi$.

(iii) $\varphi|(\psi|\sigma)\sim_{Asso}(\varphi|\psi)|\sigma$.
\end{Cor}
It follows that {\em when confined to truth in associative structures,} one can  drop parentheses from $(\varphi|\psi)|\sigma$ and write simply $\varphi|\psi|\sigma$ (as with the case of $\wedge$ and $\vee$ of in classical PL), and more generally $\varphi_1|\cdots|\varphi_n$ for any sentences $\varphi_i$ of $L_s$.

Moreover, in view of Theorem \ref{T:LO} and Corollary \ref{C:allprop}, when the choice function $f$ is associative, then $f=\min_{<}$ for a total ordering $<$ of $Sen(L)$.  Namely, the  following is proved by an easy induction:

\begin{Cor} \label{C:genarlpick}
Let $\langle M,f\rangle$ be associative, with $f=\min_{<}$ for a total ordering $<$ of $Sen(L)$. Then for every $n\in\N$ and any
 $\{\varphi_1,\ldots,\varphi_n\}\subset L_s$,
$$\langle M,f\rangle\models_s \varphi_1|\cdots|\varphi_n \ \mbox{iff} \ M\models f(\overline{f}(\varphi_1),\ldots,\overline{f}(\varphi_n)) \ \mbox{iff} \ M\models \min_<(\overline{f}(\varphi_1),\ldots,\overline{f}(\varphi_n)),$$
where $f(\sigma_1,\ldots,\sigma_n)$ abbreviates $f(\{\sigma_1,\ldots,\sigma_n\})$. In particular, for classical sentences $\alpha_1,\ldots,\alpha_n$,
$$\langle M,f\rangle\models_s \alpha_1|\cdots|\alpha_n \  \ \mbox{iff} \ M\models\min_<(\alpha_1,\ldots,\alpha_n).$$
\end{Cor}

\subsection{Regularity}

For every $\varphi\in Sen(L)$, let $Sub(\varphi)$ denote the set of sub-sentences of $\varphi$. Given   $\varphi$,  $\sigma\in Sub(\varphi)$ and any $\sigma'$, let $\varphi[\sigma'/\sigma]$ denote the result of replacing $\sigma$ by $\sigma'$ throughout $\varphi$.

\begin{Def} \label{D:Subclosed}
{\em For $X\subseteq {\cal F}$, $\sim_X$  is said to be} logically closed {\em if for all $\varphi$, $\sigma\in Sub(\varphi)$ and $\sigma'$,
$$\sigma\sim_X\sigma' \ \Rightarrow \ \varphi\sim_X\varphi[\sigma'/\sigma].$$
}
\end{Def}

Classical logical equivalence $\sim$ is logically closed of course, but $\sim_s$ and  $\sim_{Asso}$ are  not in general. The question is what further condition on $X$ is required  in order for $\sim_X$ to be  logically closed. This is the condition of {\em regularity} introduced below.

Regularity is a condition independent from  associativity, yet  compatible with it. So it is reasonable to introduce it   independently from associativity.

\begin{Def} \label{D:regular}
{\em A choice function $f$ for $L$ is said to be} regular {\em if for all  $\alpha$, $\alpha'$, $\beta$,
$$\alpha\sim\alpha' \ \Rightarrow \ \ f(\alpha,\beta)\sim f(\alpha',\beta).$$
}
\end{Def}
The following properties are immediate consequences of the definition.

\begin{Fac} \label{F:regular}
Let $f$ be regular. Then for all $\alpha$, $\alpha'$, $\beta$, $\beta'$:

(i) If   $\alpha\sim\alpha'\not\sim\beta\sim\beta'$ and $f(\alpha,\beta)=\alpha$ then $f(\alpha',\beta')=\alpha'$, while if   $f(\alpha,\beta)=\beta$ then $f(\alpha',\beta')=\beta'$.

(ii) If $\alpha\sim\alpha'\sim\beta\sim\beta'$, $f(\alpha,\beta)$ and $f(\alpha',\beta')$ can be
any element of the sets $\{\alpha,\beta\}$, $\{\alpha',\beta'\}$, respectively.
\end{Fac}

Let
$$Reg=\{f\in {\cal F}: \ \mbox{$f$ is regular}\}$$
denote  the set of regular choice functions, and let $\models_{Reg}$, $\sim_{Reg}$ abbreviate the relations $\models_{Reg}$, $\sim_{Reg}$, respectively.  Regularity not only guarantees that $\sim_{Reg}$ is logically closed, but also the  converse is true.

\begin{Thm} \label{T:regularity}
 $\sim_X$  is logically closed  if and only if  $X\subseteq Reg$.
\end{Thm}

{\em Proof.} ``$\Leftarrow$'' : We assume  $X\subseteq Reg$ and show that  $\sim_X$  is logically closed. For $\sigma\in Sub(\varphi)$, $\varphi[\sigma'/\sigma]$ is defined by induction on the length of $\varphi$ as usual, so we prove
\begin{equation} \label{E:sub}
\sigma\sim_X\sigma' \ \Rightarrow \ \varphi[\sigma'/\sigma]\sim_X\varphi
\end{equation}
along the steps of the definition of $\varphi[\sigma'/\sigma]$. Actually regularity is needed only for the treatment of case $\varphi=\varphi_1|\varphi_2$, since $\overline{f}$ commutes with standard connectives, so let us treat this step of the induction only. That is,
let $\varphi=\varphi_1|\varphi_2$, so
$$\varphi[\sigma'/\sigma]=\varphi_1[\sigma'/\sigma]|\varphi_2[\sigma'/\sigma],$$
and assume  the claim holds for $\varphi_1$, $\varphi_2$. Let us set for readability, $\varphi'=\varphi[\sigma'/\sigma]$, $\varphi_1'=\varphi_1[\sigma'/\sigma]$, $\varphi_2'=\varphi_2[\sigma'/\sigma]$. Then $\varphi'=\varphi_1'|\varphi_2'$, and by the induction assumption $\varphi_1'\sim_X\varphi_1$, $\varphi_2'\sim_X\varphi_2$. We have to show that $\varphi_1'|\varphi_2'\sim_X\varphi_1|\varphi_2$. Pick   $f\in X$. By our assumption $f\in Reg$. Our induction assumptions become
\begin{equation} \label{E:one|}
\overline{f}(\varphi_1')\sim \overline{f}(\varphi_1) \ \mbox{and} \ \overline{f}(\varphi_2')\sim \overline{f}(\varphi_2),
\end{equation}
and it suffices to show that$\overline{f}(\varphi_1'|\varphi_2')\sim \overline{f}(\varphi_1|\varphi_2)$, or equivalently,
\begin{equation} \label{E:two|}
f(\overline{f}(\varphi_1'),\overline{f}(\varphi_2'))\sim f(\overline{f}(\varphi_1),\overline{f}(\varphi_2)).
\end{equation}
But since $f\in Reg$, in view of Fact \ref{F:regular}, (\ref{E:two|}) follows immediately from (\ref{E:one|}). This completes the proof of  direction $\Leftarrow$.

``$\Rightarrow$'': Suppose  $X\not\subseteq Reg$. We have  to show that $\sim_X$ is not logically closed.  Pick $f\in X-Reg$. Since $f$  is not regular, there exist $\alpha$, $\alpha'$ and  $\beta$ in $Sen(L)$  such that $\alpha\sim\alpha'$ and  $f(\alpha,\beta)\not\sim f(\alpha',\beta)$. In particular this implies that $\alpha\not\sim \beta$. Moreover, either $f(\alpha,\beta)=\alpha$ and $f(\alpha',\beta)=\beta$, or $f(\alpha,\beta)=\beta$ and $f(\alpha',\beta)=\alpha'$.

Assume $f(\alpha,\beta)=\alpha$ and $f(\alpha',\beta)=\beta$, the other case being treated similarly. Since   $\alpha\not\sim \beta$ there exists a truth assignment $M$ for $L$ such that $M\models\alpha\wedge\neg\beta$ or $M\models\neg\alpha\wedge\beta$. Without loss of generality assume $M\models\alpha\wedge\neg\beta$.  Then  $M\models f(\alpha,\beta)$ which means $\langle M,f\rangle\models\alpha|\beta$. On the other hand,  since $f(\alpha',\beta)=\beta$ and  $M\models\neg\beta$, we have  $M\models\neg f(\alpha',\beta)$, that is, $M\not\models f(\alpha',\beta)$, which means  $\langle M,f\rangle\not\models\alpha'|\beta$. Therefore for some $M$ and some $f\in X$,
$\langle M,f\rangle\models\alpha|\beta$ and  $\langle M,f\rangle\not\models\alpha'|\beta$. Thus $\alpha'|\beta\not\sim_X\alpha|\beta$, while $\alpha'\sim \alpha$, and hence $\alpha'\sim_X\alpha$. It follows  that $\sim_X$ is not logically closed.  \telos

\vskip 0.2in

In general if  $X\subseteq Y\subseteq {\cal F}$ and $\sim_X$ is logically closed, it doesn't  seem likely that one can infer that $\sim_Y$ is logically closed (or vice-versa). Yet we have the following generalization of  \ref{T:regularity}.

\begin{Cor} \label{C:compare}
If $X\subseteq Reg\subseteq Y\subseteq {\cal F}$, then for all $\varphi$, $\sigma\in Sub(\varphi)$ and  $\sigma'$,
$$\sigma\sim_Y\sigma' \ \Rightarrow \ \varphi[\sigma'/\sigma]\sim_X\varphi.$$
\end{Cor}

{\em Proof.} By Fact \ref{F:tautolol} (iii), $X\subseteq Y$ implies $\sim_Y\subseteq \sim_X$. So given $X\subseteq Reg\subseteq Y$, we have $\sim_Y\subseteq\sim_{Reg}\subseteq\sim_X$. Thus if $\sigma\sim_Y\sigma'$, then  $\sigma\sim_{Reg}\sigma'$. By theorem \ref{T:regularity}, this implies $\varphi[\sigma'/\sigma]\sim_{Reg}\varphi$, and therefore $\varphi[\sigma'/\sigma]\sim_X\varphi$. \telos

\vskip 0.2in

Using the axiom of choice one can easily construct  regular choice functions for  $L$.  For every $\alpha\in Sen(L)$, let $[\alpha]$ denote the $\sim$-equivalence class of $\alpha$, i.e.,
$$[\alpha]=\{\beta:\beta\sim\alpha\}.$$

\begin{Prop} \label{P:ac}
{\rm (AC)} There exist  regular choice functions for $L$.
\end{Prop}

{\em Proof.}  Using $AC$, pick a representative $\xi_\alpha$ from each equivalence class $[\alpha]$ and  let $D=\{\xi_\alpha:\alpha\in Sen(L)\}$. Every $\alpha\in Sen(L)$ is logically equivalent to  $\xi_\alpha\in D$.  Let $f_0:[D]^2\rightarrow D$ be an arbitrary choice function for all pairs of elements of $D$. Then $f_0$ extends  to a regular choice function $f$ for $L$, defined as follows:

$f(\alpha,\beta)=\alpha$, if $\alpha\not\sim \beta$ and $f_0(\xi_\alpha,\xi_\beta)=\xi_\alpha$. If $\alpha\sim\beta$, we define $f(\alpha,\beta)$ arbitrarily (to be precise, by setting $f(\alpha,\beta)=g(\alpha,\beta)$, where $g$ is  a choice function  on all pairs $\{\alpha,\beta\}$ of sentences such that $\alpha\sim \beta$).
\telos

\vskip 0.2in

Next we come  to {\em associative  regular} choice functions.

\begin{Def} \label{D:regularord}
{\em A total ordering $<$ of $Sen(L)$ is said to be} regular {\em if for all  $\alpha,\beta$,
$$\alpha\not\sim \beta \ \& \ \alpha<\beta \ \Rightarrow \ [\alpha]<[\beta],$$
(where  $[\alpha]<[\beta]$  means  that for all $\alpha'\in [\alpha]$ and $\beta'\in[\beta]$, $\alpha'<\beta'$).
}
\end{Def}
The following is an immediate consequence of the  preceding definitions.

\begin{Fac} \label{F:assregular}
Let $<$ be a total ordering of $Sen(L)$. Then $<$ is regular if and only if the corresponding associative choice function $f=\min_<$ is regular.
\end{Fac}

Thus the following simple construction of  regular total orderings for  $Sen(L)$ supplements Proposition \ref{P:ac} above.

\begin{Prop} \label{P:plenty}
{\rm (AC)} (i) There  exist regular total orderings of $Sen(L)$.

(ii) Moreover, for any set $A\subset Sen(L)$ of pairwise inequivalent sentences, and any partial ordering $R$ of $A$, there is a regular total ordering $<$ of $Sen(L)$ such that  $R\subseteq <$.
\end{Prop}

{\em Proof.} (i)  Let $Sen(L)/\!\!\sim$ be the set of equivalence classes $[\alpha]$, $\alpha\in Sen(L)$. For each $[\alpha]\in Sen(L)/\!\!\sim$ pick by $AC$ a total ordering $<_{[\alpha]}$ of $[\alpha]$. Pick also  a total ordering $<_1$ of $Sen(L)/\!\!\sim$. These orderings generate  a regular total ordering $<$ of $Sen(L)$ defined by:

$\alpha<\beta$ if and only if $\alpha\not\sim\beta$ and $[\alpha]<_1[\beta]$, or $\alpha\sim\beta$ and $\alpha<_{[\alpha]} \beta$.

\noindent (ii) Since the  elements of $A$ are pairwise inequivalent, we can think of $A$ as a subset of $Sen(L)/\sim$. Since $R$ is already a partial  ordering of $A$,  it suffices to pick (by the help of $AC$) the total ordering $<_1$ of the preceding case so that $R\subset <_1$.    \telos

\vskip 0.2in

Both associativity and regularity are indispensable for a reasonable notion of truth $\models_X$ that captures the behavior of $|$. This is because on the one  hand  without associativity  one would have to face unmanageable  complexity caused by incomparable sentences  of the form $((\alpha|\beta)|\gamma)|\delta$, $(\alpha|\beta)|(\gamma|\delta)$, $\alpha|((\beta|\gamma)|\delta)$, etc.
On the other hand regularity entails  logical closeness,  without which one  cannot establish  even that the sentences, for example,  $\alpha|\beta$ and $\alpha|\neg\neg\beta$ are essentially identical.  Thus a natural class of choice functions to work with is
$$Reg^*=Reg \cap Asso.$$
We abbreviate the corresponding semantic notions $\models_{Reg^*}$, $\sim_{Reg^*}$, by $\models_{Reg^*}$ and  $\sim_{Reg^*}$, respectively.

\vskip 0.2in

Note that in view of regularity (and only in view of that) we can write, for example, $\varphi|\top$ and $\varphi|\bot$, where $\top$ and $\bot$ denote the classes of classical tautologies and contradictions, respectively.

The next question is how the standard connectives act on $|$ and vice-versa. Specifically we shall examine  whether:

(a) $\neg$ can commute with $|$,

(b) $\wedge$ and $\vee$ can distribute over $|$,

(c) $|$ can distribute over $\wedge$ and $\vee$.

\noindent We shall see that all three questions are answered in the negative with respect to the truth relations $\models_{Reg^*}$.

Concerning the first question one  can construct  a choice function $f$ such that for every truth assignment $M$  and any sentences $\varphi$, $\psi$, $$\langle M,f\rangle \models_s\neg(\varphi|\psi)\leftrightarrow \neg\varphi|\neg\psi.$$
For that it suffices to define $f$ so that $\overline{f}(\neg\varphi|\neg\psi)=\overline{f}(\neg(\varphi|\psi))$,
or equivalently
$$f(\neg\overline{f}(\varphi),\neg\overline{f}(\psi))=\neg f(\overline{f}(\varphi), \overline{f}(\psi)).$$
This can be done by defining $f(\alpha,\beta)$ by induction along an ordering of the pairs $\langle r(\alpha),r(\beta)\rangle$, where $r(\alpha)$ is the usual rank of $\alpha$.

Nevertheless,  such an $f$ supporting $\neg(\varphi|\psi)\leftrightarrow \neg\varphi|\neg\psi$ would  serve just as a counterexample or a curiosity,  and could by no means characterize a natural class of functions. Specifically it is easily seen that such an $f$ cannot be regular.

\begin{Fac} \label{F:curious}
If $f$ is  regular then for every $M$ there are $\varphi$, $\psi$ such that $\langle M,f\rangle\not\models \neg(\varphi|\psi)\leftrightarrow \neg\varphi|\neg\psi$. Thus for regular $f$, the scheme $\neg(\varphi|\psi)\leftrightarrow \neg\varphi|\neg\psi$ is always false in $\langle M,f\rangle$.
\end{Fac}

{\em Proof.} Let $f$ be regular. Then for every $\alpha$,  $f(\neg\alpha,\neg\neg\alpha)\sim f(\neg\alpha,\alpha)=f(\alpha,\neg\alpha)$, therefore $\neg f(\alpha,\neg\alpha)\leftrightarrow f(\neg\alpha,\neg\neg\alpha)$ is a contradiction. Hence for every $M$, $M\not\models \neg f(\alpha,\neg\alpha)\leftrightarrow f(\neg\alpha,\neg\neg\alpha)$, which means that $M\not\models \overline{f}(\neg (\alpha|\neg\alpha)\leftrightarrow \neg\alpha|\neg\neg\alpha)$, or
$\langle M,f\rangle \not\models \neg (\alpha|\neg\alpha)\leftrightarrow \neg\alpha|\neg\neg\alpha$. \telos

\vskip 0.2in

Concerning question (b) above the answer is negative with  respect to the semantics $\models_X$ for any $X\subseteq Reg^*$. Let us give some definitions with the purpose to prove later that they are void.

\begin{Def} \label{D:monotonic}
{\em Let $<$ be a regular total ordering of $Sen(L)$.  $<$ is said to be:

(a)} $\wedge$-monotonic, {\em if  for all $\alpha,\beta,\gamma\in Sen(L)$
such that  $\alpha\wedge\gamma\not\sim\beta\wedge\gamma$,
$$\alpha<\beta \  \Leftrightarrow \alpha\wedge\gamma<\beta\wedge\gamma.$$

(b)} $\vee$-monotonic, {\em if   for all $\alpha,\beta,\gamma\in Sen(L)$
such that $\alpha\vee\gamma\not\sim\beta\vee\gamma$,
$$\alpha<\beta \ \Leftrightarrow \alpha\vee\gamma<\beta\vee\gamma.$$

Accordingly, a choice function $f\in Reg^*$ is said to be} $\wedge$-monotonic {\em (resp.} $\vee$-monotonic{\em) if $f=\min_<$ and $<$ is $\wedge$-monotonic (resp. $\vee$-monotonic).
}
\end{Def}

\begin{Lem} \label{L:ifmonoton}
(i) If $<$ is $\wedge$-monotonic, then  $$\min(\alpha\wedge\gamma,\beta\wedge\gamma)\sim\gamma\wedge\min(\alpha,\beta).$$

(ii) If $<$ is $\vee$-monotonic, then
$$\min(\alpha\vee\gamma,\beta\vee\gamma)\sim\gamma\vee\min(\alpha,\beta).$$

\end{Lem}

{\em Proof.} (i) If  $\alpha\wedge\gamma\sim \beta\wedge\gamma$ then  obviously        $\min(\alpha\wedge\gamma,\beta\wedge\gamma)\sim\gamma\wedge \ \min(\alpha,\beta)$.
So assume  $\alpha\wedge\gamma\not\sim\beta\wedge\gamma$. Then also $\alpha\not\sim\beta$. Without loss of generality  suppose $\alpha<\beta$, so $\min(\alpha,\beta)=\alpha$. By $\wedge$-monotonicity, $\alpha\wedge\gamma<\beta\wedge\gamma$, so
$\min(\alpha\wedge\gamma,\beta\wedge\gamma)=\alpha\wedge\gamma=
\gamma\wedge\min(\alpha,\gamma)$.

(ii) Similar.   \telos

It is easy to  give syntactic characterizations of $\wedge$- and $\vee$-monotonicity. The proof of the following is left to the reader.

\begin{Lem} \label{L:montruth}
Let $f\in Reg^*$. Then:

(i) $f$ is $\wedge$-monotonic if and only if  for all $M$ and all $\varphi,\psi,\sigma\in Sen(L_s)$:
$$\langle M,f\rangle\models\varphi\wedge(\psi|\sigma)\leftrightarrow (\varphi\wedge \psi)|(\varphi\wedge\sigma).$$

(ii)  $f$ is $\vee$-monotonic if and only if for all $M$ and all $\varphi,\psi,\sigma\in Sen(L_s)$:
$$\langle M,f\rangle\models\varphi\vee(\psi|\sigma)\leftrightarrow (\varphi\vee \psi)|(\varphi\vee\sigma).$$

\end{Lem}

It follows from the previous Lemma that $\wedge$- and $\vee$-monotonicity are exactly  the conditions under which $\wedge$ and $\vee$, respectively, distribute over $|$. However we can easily see by a counterexample that there are no $\wedge$-monotonic or $\vee$-monotonic regular functions (or orderings).

\begin{Prop} \label{F:nomonotonic}
There is no regular total ordering $<$ of $Sen(L)$ which is $\wedge$-monotonic or $\vee$-monotonic. Consequently there is no $X\subseteq Reg^*$ such that the schemes $\varphi\wedge(\psi|\sigma)\leftrightarrow (\varphi\wedge \psi)|(\varphi\wedge\sigma)$ and $\varphi\vee(\psi|\sigma)\leftrightarrow (\varphi\vee \psi)|(\varphi\vee\sigma)$, are not $\models_X$-tautologies.
\end{Prop}

{\em Proof.} Suppose there is a regular total ordering $<$ of $Sen(L)$ which is $\wedge$-monotonic. Let $p$, $q$ $r$ be atomic sentences such that $p<q<r$. Consider the formula $\alpha=p\wedge r\wedge \neg q$. Then by $\wedge$- monotonicity $p<q$ implies $p\wedge r\wedge \neg q<q\wedge r\wedge \neg q$, or by regularity $\alpha<\bot$. For the same reason $q<r$ implies $p\wedge q\wedge \neg q<p\wedge r\wedge \neg q$, or $\bot<\alpha$, a contradiction. Working with $\beta=p\vee r\vee \neg q$ we similarly show that  is no regular total ordering $<$  which is $\vee$-monotonic. \telos

\vskip 0.2in

Having settled the question about the  distributivity of $\wedge$ and $\vee$ over $|$, we come to the converse question, whether $|$ can distribute over $\wedge$ and/or $\vee$ for some class $X$ of choice functions such that $X\subseteq Reg^*$. The answer is ``no''  again with respect to $Reg^*$. Namely:

\begin{Prop} \label{F:supdistrib}
There is no  regular total ordering $<$ of $Sen(L)$ such that if $f=\min_<$, then  for every $M$, $\langle M,f\rangle$ satisfies the scheme
$$(*) \ \ \ \varphi|(\psi\wedge \sigma)\leftrightarrow (\varphi|\psi)\wedge (\varphi|\sigma).$$
Consequently there is no $X\subseteq Reg^*$ such that (*) is a $\models_X$-tautology. Similarly for the dual scheme
$$(**) \ \ \ \varphi|(\psi\vee \sigma)\leftrightarrow (\varphi|\psi)\vee (\varphi|\sigma).$$
\end{Prop}

{\em Proof.} Towards reaching a  contradiction assume that  there is a regular total ordering $<$ of $Sen(L)$ such that if $f=\min_<$, then   $(*)$ is true in all models $\langle M,f\rangle$. Fix some atomic sentence $p$ of $L$. By regularity we have either $p<\bot$ or $\bot<p$. We examine below some consequences of each of these  cases.

(i) Let $p<\bot$. Pick some $q\neq p$ and an $M$ such that $M\models p\wedge q$. Then
$$\langle M,f\rangle\models_s p|(q\wedge \neg q)\leftrightarrow (p|q)\wedge (p|\neg q),$$
or
\begin{equation} \label{E:basic1}
\langle M,f\rangle\models_s p|\bot\leftrightarrow (p|q)\wedge (p|\neg q).
\end{equation}
Since $p<\bot$ and $M\models p$, the left-hand side of the equivalence in (\ref{E:basic1}) is true in $\langle M,f\rangle$. Thus so is the right-hand side of the equivalence. Since $M\models p\wedge q$, the conjunct $p|q$ is true, while the truth of the  conjunct $p|\neg q$ necessarily implies  $p<\neg q$, since $M\models q$. Then pick $N$ such that $N\models p\wedge \neg q$. We have also
\begin{equation} \label{E:basic2}
\langle N,f\rangle\models_s p|\bot\leftrightarrow (p|q)\wedge (p|\neg q).
\end{equation}
Again the left-hand side of the equivalence in (\ref{E:basic2}) is true in $\langle N,f\rangle$.  So the right-hand side  is true too. Since $N\models p\wedge \neg q$, the conjunct $p|\neg q$ holds. In order for the conjunct $p|q$ to hold too we must have $p< q$, since $N\models \neg q$.  Summing up the above two facts we conclude that if the letters $p$, $q$ range over atomic sentences, then
\begin{equation} \label{E:fact1}
(\forall p\neq q)(p<\bot \ \Rightarrow p<q \ \& \ p<\neg q).
\end{equation}

(ii) Let now $\bot<p$. Pick again some $q\neq p$ and an $M$ such that $M\models p\wedge q$. Then (\ref{E:basic1}) holds again, but now the left-hand side of the equivalence in (\ref{E:basic1}) is false $\langle M,f\rangle$. Thus so is the right-hand side, which, since $M\models p$, necessarily implies $\neg q<p$. Then  pick $N$ such that $N\models p\wedge \neg q$.  (\ref{E:basic2}) holds again with the left-hand side of the equivalence being false. The right-hand side is false too and this holds only if $q<p$, since $N\models\neg q$. Therefore from these two facts we conclude that
\begin{equation} \label{E:fact2}
(\forall p\neq q)(\bot<p \ \Rightarrow q<p \ \& \ \neg q<p).
\end{equation}
Now  since there are at least three distinct atoms $p,q,r$ and $p,q,r,\bot$ are linearly ordered by $<$, then at least two of them lie on the left of $\bot$, or on the right of $\bot$. That is, there are $p,q$ such that $p,q<\bot$ or $\bot<p,q$. If $p,q<\bot$, (\ref{E:fact1}) implies that $p<q$, $p<\neg q$, $q<p$ and $q<\neg p$, a contradiction. If $\bot<p,q$, then  (\ref{E:fact2}) implies that $q<p$, $\neg q<p$, $p<q$ and $\neg p<q$, a contradiction again. This completes the proof that (*) cannot be a $\models_X$-tautology for any $X\subseteq Reg^*$. Concerning the scheme (**) we consider the instances
$$p|(q\vee\neg q)\leftrightarrow (p|q)\vee (q|\neg q),$$
i.e.,
$$p|\top\leftrightarrow (p|q)\vee (q|\neg q),$$
for atomic sentences $p,q$, and argue analogously as before, by examining the cases $p<\top$ and $\top<p$. \telos

\subsection{$\neg$-decreasingness}
There is still the question of how $\neg$ behaves with respect to $|$. As we saw in Fact \ref{F:curious}, $\neg$ cannot commute with $|$ in models $\langle M,f\rangle$ with regular $f$. Equivalently, if $f\in Reg^*$ and $f=\min_<$,  $\neg$ cannot be ``increasing'', that is, cannot   satisfy $\alpha<\beta\Leftrightarrow \neg\alpha<\neg\beta$, for all $\alpha$, $\beta$. However it can be ``decreasing'', and this turns out to be a useful property.

\begin{Def} \label{D:negdec}
{\em $<$ is said to be} $\neg$-decreasing, {\em if  for all $\alpha,\beta\in Sen(L)$ such that $\alpha\not\sim\beta$,
$$\alpha<\beta \Leftrightarrow \neg\beta<\neg\alpha.$$
Accordingly a choice function $f\in Reg^*$ is said to be} $\neg$-decreasing {\em if $f=\min_<$ and $<$ is $\neg$-decreasing.}
\end{Def}

\begin{Lem} \label{L:decrease}
$<$ is $\neg$-decreasing if and only if for all $\alpha\not\sim\beta$,  $$\neg\min(\neg\alpha,\neg\beta)\sim\max(\alpha,\beta).$$
\end{Lem}

{\em Proof.} Let $<$ be $\neg$-decreasing. Then for any $\alpha\not\sim\beta$, $\alpha<\beta \Leftrightarrow \neg\beta<\neg\alpha$, so, $\min(\neg\alpha,\neg\beta)=\neg\max(\alpha,\beta)$, hence
$\neg\min(\neg\alpha,\neg\beta)\sim\max(\alpha,\beta)$.

Conversely, suppose $<$ is not  $\neg$-decreasing. Then there are $\alpha$, $\beta$ such that $\alpha\not\sim\beta$, $\alpha<\beta$ and $\neg\alpha<\neg\beta$. But then $\neg\min(\neg\alpha,\neg\beta)=\neg\neg\alpha\not\sim\beta=\max(\alpha,\beta)$.\telos

\vskip 0.2in

We can also give a syntactic characterization of $\neg$-decreasingness.

\begin{Thm} \label{T:syncar}
Let $f\in Reg^*$. Then $f$ is $\neg$-decreasing if and only if for every $M$ and any $\varphi$, $\psi$,
$$\langle M,f\rangle\models_s\varphi\wedge\neg\psi\rightarrow (\varphi|\psi\leftrightarrow \neg\varphi|\neg\psi).$$
\end{Thm}

{\em Proof.} ``$\Rightarrow$'': Let $f$ be $\neg$-decreasing and $f=\min_<$. Let $M$ and $\varphi$, $\psi$ such that $\langle M,f\rangle\models_s\varphi\wedge\neg\psi$, that is, $M\models \overline{f}(\varphi)\wedge\neg\overline{f}(\psi)$.
It suffices to show that
$\langle M,f\rangle\models_s(\varphi|\psi\leftrightarrow \neg\varphi|\neg\psi)$, or equivalently,
$$M\models f(\overline{f}(\varphi),\overline{f}(\psi))\leftrightarrow f(\neg\overline{f}(\varphi),\neg\overline{f}(\psi)).$$
If we set $\overline{f}(\varphi)=\alpha$ and $\overline{f}(\psi)=\beta$, the above amount to assuming  that $M\models\alpha\wedge\neg\beta$ and concluding  that $M\models f(\alpha,\beta)\leftrightarrow f(\neg\alpha,\neg\beta)$, or $$M\models \min(\alpha,\beta)\leftrightarrow \min(\neg\alpha,\neg\beta).$$ But since $M\models\alpha\wedge\neg\beta$, $M\models \min(\alpha,\beta)$ implies
$\min(\alpha,\beta)=\alpha$. Then, since $<$ is $\neg$-decreasing, $\min(\neg\alpha,\neg\beta)=\neg\beta$. Therefore $M\models \min(\neg\alpha,\neg\beta)$. So $M\models \min(\alpha,\beta)\rightarrow \min(\neg\alpha,\neg\beta)$. The converse is similar.

``$\Leftarrow$'': Let $f$ be  non-$\neg$-decreasing, with $f=\min_<$,  and let  $\alpha\not\sim\beta$ such that    $\alpha<\beta$ and $\neg\alpha<\neg\beta$.  Without loss of generality there is $M$ such that $M\models\neg\alpha\wedge\beta$. Then $M\not\models \min(\alpha,\beta)$, thus $\langle M,f\rangle\not\models_s \alpha|\beta$, while $M\models \min(\neg\alpha,\neg\beta)$, or $\langle M,f\rangle\models_s \neg\alpha|\neg\beta$. So $\langle M,f\rangle\not\models (\alpha|\beta\leftrightarrow \neg\alpha|\neg\beta)$, and therefore
$$\langle M,f\rangle\not\models \neg\alpha\wedge\beta\rightarrow (\alpha|\beta\leftrightarrow \neg\alpha|\neg\beta).$$
Therefore $\langle M,f\rangle$ does not satisfy the scheme $\varphi\wedge\neg\psi\rightarrow (\varphi|\psi\leftrightarrow \neg\varphi|\neg\psi)$. \telos

\vskip 0.2in

Next  let us make sure that $\neg$-decreasing total orderings exist.

\begin{Thm} \label{T:existdecreas}
There exist  regular $\neg$-decreasing total orderings of $Sen(L)$, and hence  regular $\neg$-decreasing choice functions for $L$.
\end{Thm}

{\em Proof.} There is a general method for constructing regular and $\neg$-decreasing total orderings of $Sen(L)$ that makes use of  the Axiom of Choice. This is the following. Let $P=\{\{[\alpha],[\neg\alpha]\}:\alpha\in Sen(L)\}$. Pick by AC a choice function $F$ for $P$, and let $A=\bigcup F``P$ and $B=Sen(L)-A$. Both $A$, $B$ are $\sim$-saturated, that is, $\alpha\in A\Rightarrow [\alpha]\subset A$, and similarly for $B$.
As in the proof of Proposition \ref{P:plenty} pick a regular total ordering $<_1$ of $A$. By the definition of $A$, $B$, clearly for every $\alpha\in Sen(L)$, $$\alpha\in A\Leftrightarrow \neg\alpha\in B,$$
so $<_1$ induces a regular total ordering $<_2$ of $B$ by setting
$$\alpha<_2\beta \Leftrightarrow \neg\beta<_1\neg\alpha.$$
Then define $<$ of $Sen(L)$ as follows: $\alpha<\beta$ if and only if :

$\alpha\in A$ and $\beta\in B$, or

$\alpha,\beta\in A$ and $\alpha<_1\beta$, or

$\alpha,\beta\in B$ and $\alpha<_2\beta$.

\noindent It is easy to verify that $<$ is a total regular and $\neg$-decreasing ordering of $Sen(L)$. \telos

\vskip 0.2in

We can further show that every regular $\neg$-decreasing total ordering is constructed by the general method of Theorem  \ref{T:existdecreas}. Let us give some definitions. A set $X\subset Sen(L)$ is said to be {\em selective} if of every pair of opposite sentences $\{\alpha,\neg\alpha\}$ $X$ contains exactly one. Recall also that $X$ is {\em $\sim$-saturated} if for every $\alpha$, $\alpha\in X\Rightarrow [\alpha]\subset X$. Note that  familiar examples of  selective and $\sim$-saturated sets are  the consistent and complete sets $\Sigma\subset Sen(L)$ (as well as their complements $Sen(L)-\Sigma$). However not every selective and $\sim$-saturated set is of this kind.
For instance the sets $A$, $B$ in the proof of \ref{T:existdecreas} are selective and $\sim$-saturated forming a partition of $Sen(L)$. Moreover $A$ is an initial segment and $B$ is a final segment of  $\langle Sen(L),<\rangle$.  We shall see that such a partition exists  for every regular $\neg$-decreasing total ordering.

\begin{Prop} \label{P:segments}
Let $<$ be  a regular $\neg$-decreasing total ordering of $Sen(L)$. Then $Sen(L)$ splits into two $\sim$-saturated  sets $I$ and $J$ which are  selective, hence
$$\alpha\in I \Leftrightarrow \neg\alpha\in J,$$
and  $I<J$, that is, $I$ is an initial and $J$ a final segment of $<$.
\end{Prop}

{\em Proof.} Let $<$ be  a regular and $\neg$-decreasing total ordering of $Sen(L)$. Let us call an initial segment of $<$ {\em weakly selective} if from  every pair $\{\alpha,\neg\alpha\}$, $I$ contains {\em at most one} element. We first claim that there are weakly selective initial segments of $<$. Observe  that if for some $\alpha$,  $(\forall \beta)(\alpha<\beta \vee \alpha<\neg\beta)$ is true, then the initial segment $\{\beta:\beta\leq\alpha\}$ is weakly selective. So if, towards a contradiction, we assume that no weakly selective initial segment exists, then
\begin{equation} \label{E:select}
\forall\alpha \exists \beta (\beta\leq \alpha \wedge \neg\beta\leq \alpha).
\end{equation}
Assume (\ref{E:select}) holds and fix some $\alpha$. Pick $\beta$ such that $\beta\leq \alpha$ and $\neg\beta\leq \alpha$. By $\neg$-decreasingness, $\neg\alpha\leq \neg\beta$. Therefore $\neg\alpha<\alpha$. Now apply (\ref{E:select}) to $\neg\alpha$ to find $\gamma$ such that $\gamma\leq \neg\alpha$ and $\neg\gamma\leq \neg\alpha$. By $\neg$-decreasingness, $\neg\neg\alpha\leq \neg\gamma$. Thus $\neg\neg\alpha\leq \neg\alpha$ and  by regularity, $\alpha\leq \neg\alpha$. But this contradicts $\neg\alpha<\alpha$.

So there exist weakly selective initial segments of $<$. Taking the union of all such initial  segments, we find a greatest weakly selective initial segment $I$. It is easy to see that $I$ is selective, i.e., from each pair $\{\alpha,\neg\alpha\}$ it contains {\em exactly one} element. Indeed, assume the contrary. Then there is $\alpha$, such that either $I<\alpha<\neg\alpha$, or $I<\neg\alpha<\alpha$. Assume the first is the case, the other being similar. But then there is $\beta$ such that $I<\{\beta,\neg\beta\}<\alpha<\neg\alpha$, because otherwise the segment $\{\gamma:\gamma\leq \alpha\}$ would be a weakly selective segment greater than $I$, contrary to  the maximality of $I$. Now $\{\beta,\neg\beta\}<\alpha<\neg\alpha$ implies that $\beta<\alpha$ and $\neg\beta<\neg\alpha$, which contradicts the $\neg$-decreasingness of $<$.

Further, let $J=\{\neg\alpha:\alpha\in I\}$. Then $J$ is a greatest selective final segment, $I<J$ and $I\cap J=\emptyset$. To show that $I$ (and hence $J$) is $\sim$-saturated, let $\alpha\in I$. Assume first that $\alpha$ is not the  greatest element of $I$, so  there is  $\beta\in I$ such that $\alpha<\beta$. By regularity, $[\alpha]<\beta$. Hence $[\alpha]\subset I$. Next assume that $\alpha$ is the greatest element of $I$. Then necessarily $\alpha$ is the greatest element of $[\alpha]$ too, otherwise $I\cup[\alpha]\supsetneq I$ and $I\cup [\alpha]$ is selective, contrary to the maximality of $I$. Thus again $I$ is $\sim$-saturated.

It remains  to show that $I\cup J=Sen(L)$. Assume $\alpha\notin I$. Since $I$ is selective, $\neg\alpha\in I$, therefore $\neg\neg\alpha\in J$. Since $J$ is $\sim$-saturated, $\alpha\in J$. Thus $I\cup J=Sen(L)$.  \telos

\vskip 0.2in

So regular $\neg$-decreasing functions constitute a natural class of choice functions stronger than $Reg^*$. Let
$$Dec=\{f\in Reg^*: f \ \mbox{is $\neg$-decreasing}\}.$$
We abbreviate the corresponding semantic notions $\models_{Dec}$ and  $\sim_{Dec}$, by $\models_{Dec}$ and  $\sim_{Dec}$, respectively.

\subsection{The dual connective}
Every binary or unary logical operation when combined with negation produces  a {\em dual} one. The dual of $|$ is
$$\varphi\circ\psi:=\neg(\neg\varphi|\neg\psi)$$
for all $\varphi,\psi\in Sen(L_s)$.

A natural  question is whether each of the  operations  $|$ and $\circ$ distributes  over its dual with respect to a  truth relation $\models_X$, that is, whether there is a class of functions $X$ such that
\begin{equation} \label{E:citcdis}
\models_X\varphi\circ (\psi|\sigma)\leftrightarrow (\varphi\circ\psi)|(\varphi\circ\sigma)
\end{equation}
 and
\begin{equation} \label{E:bardis}
\models_X\varphi| (\psi\circ\sigma)\leftrightarrow (\varphi|\psi)\circ(\varphi|\sigma)
\end{equation}
are $X$-tautologies. (\ref{E:citcdis}), (\ref{E:bardis}) are dual and  equivalent to each other,  since taking the negations of both sides of (\ref{E:citcdis}) one obtains (\ref{E:bardis}), and vice-versa.

\begin{Prop} \label{P:notdistr}
There exist $\alpha$, $\beta$, $\gamma$, $M$ and $f\in Reg^*$ ($f$ non-$\neg$-decreasing) such that
$$\langle M,f\rangle\not\models_s
\alpha\circ (\beta|\gamma)\leftrightarrow(\alpha\circ\beta)|(\alpha\circ\gamma).$$
\end{Prop}

{\em Proof.} Pick  $\alpha$, $\beta$, $\gamma$  such that $\alpha\not\models\gamma$, $\gamma\not\sim\neg\alpha$, and  $M$ such that $M\models\alpha\wedge\neg\gamma$. Then we can easily find  a regular total ordering of $Sen(L)$ such that $\neg\alpha<\neg\beta$, $\neg\gamma<\neg\alpha$, $\gamma<\alpha$ and $\beta<\gamma$. Let $f=\min_<=\min$.  By definition, $\alpha\circ (\beta|\gamma)=\neg(\neg\alpha|\neg(\beta|\gamma))$, and $(\alpha\circ\beta)|(\alpha\circ\gamma)=
\neg(\neg\alpha|\neg\beta)|\neg(\neg\alpha|\neg\gamma))$. Therefore
$$\overline{f}(\alpha\circ (\beta|\gamma))=\overline{f}(\neg(\neg\alpha|\neg(\beta|\gamma))
=\neg\overline{f}(\neg\alpha|\neg(\beta|\gamma))=\neg f(\neg\alpha,\neg f(\beta,\gamma))=$$
$$\neg\min(\neg\alpha,\neg\min(\beta,\gamma))=
\neg\min(\neg\alpha,\neg\beta)=\neg\neg\alpha.$$
On the other hand,
$$\overline{f}((\alpha\circ\beta)|(\alpha\circ\gamma))=
\overline{f}(\neg(\neg\alpha|\neg\beta)|\neg(\neg\alpha|\neg\gamma))=
f(\neg f(\neg\alpha,\neg\beta),\neg f(\neg\alpha,\neg\gamma))=$$
$$\min(\neg \min(\neg\alpha,\neg\beta),\neg\min(\neg\alpha,\neg\gamma))=
\min(\neg\neg\alpha,\neg\neg\gamma)=\neg\neg\gamma,$$
where the last equation is due to the fact that $\min(\alpha,\gamma)=\gamma$ and $<$ is regular. Thus $M\models\overline{f}(\alpha\circ (\beta|\gamma))$ and $M\not\models\overline{f}((\alpha\circ\beta)|(\alpha\circ\gamma))$. Therefore
$\langle M,f\rangle\models_s
\alpha\circ (\beta|\gamma)$ and $\langle M,f\rangle\not\models_s(\alpha\circ\beta)|(\alpha\circ\gamma)$. \telos

\vskip 0.2in

Note that in the preceding counterexample we have  $\gamma<\alpha$ and $\neg\gamma<\neg\alpha$, so the  ordering $<$  is not $\neg$-decreasing. We see next that  if  $f$ is $\neg$-decreasing, then in $\langle M,f\rangle$ $|$ and $\circ$ do distribute over each other.

\begin{Prop} \label{P:distributive}
If $f\in Reg^*$ is $\neg$-decreasing, then for all  $M$, $\varphi$, $\psi$, $\sigma$,
$$\langle M,f\rangle\models_s\varphi\circ(\psi|\sigma)\leftrightarrow
(\varphi\circ\psi)|(\varphi\circ\sigma),$$
and
$$\langle M,f\rangle\models_s\varphi|(\psi\circ\sigma)\leftrightarrow
(\varphi|\psi)\circ(\varphi|\sigma).$$
\end{Prop}

{\em Proof.} The above equivalences are dual to each other, so it suffices to show the first of them. Specifically it suffices to prove that  if $f\in Reg^*$ and $f$ is $\neg$-decreasing, then
$$\overline{f}(\varphi\circ(\psi|\sigma))\sim
\overline{f}((\varphi\circ\psi)|(\varphi\circ\sigma)).$$
Fix such an $f$ and let $<$ be the  regular,  $\neg$-decreasing  ordering such that $f=\min_<=\min$. If we set  $\overline{f}(\varphi)=\alpha$, $\overline{f}(\psi)=\beta$, $\overline{f}(\sigma)=\gamma$,  express $\circ$ in terms of $|$ and replace  $f$ with $<$,
the above equivalence is written:
\begin{equation} \label{E:minmin}
\neg\min(\neg\alpha,\neg\min(\beta,\gamma))\sim
\min(\neg\min(\neg\alpha,\neg\beta),\neg\min(\neg\alpha,\neg\gamma)).
\end{equation}
If $\alpha\sim\beta\sim\gamma$, obviously (\ref{E:minmin}) is true.
Assume  $\alpha\sim\beta$ and $\alpha\not\sim\gamma$. Then,  by regularity, (\ref{E:minmin}) becomes
$$\neg\min(\neg\alpha,\neg\min(\alpha,\gamma))\sim
\min(\alpha,\neg\min(\neg\alpha,\neg\gamma)).$$
To verify it we consider the cases  $\alpha<\gamma$ and $\gamma<\alpha$. E.g. let $\alpha<\gamma$. By $\neg$-decreasingness, $\neg\gamma<\neg\alpha$, so both sides of the above relation are $\sim$ to $\alpha$. Similarly if $\gamma<\alpha$.

So it remains to prove (\ref{E:minmin}) when $\alpha\not\sim\beta$ and $\alpha\not\sim\gamma$. Then, by Lemma \ref{L:decrease}, $\neg\min(\neg\alpha,\neg\beta)\sim\max(\alpha,\beta)$,
so  (\ref{E:minmin}) is written
\begin{equation} \label{E:minmax}
\max(\alpha,\min(\beta,\gamma))\sim
\min(\max(\alpha,\beta),\max(\alpha,\gamma)).
\end{equation}
We don't know if  there is some more elegant direct (that is, not-by-cases) proof of  (\ref{E:minmax}). So we verify it by cases.

{\em Case 1.} Assume  $\alpha\leq \min(\beta,\gamma)$. Then  $\max(\alpha,\min(\beta,\gamma))=\min(\beta,\gamma)$. Besides $\alpha\leq \min(\beta,\gamma)$ implies $\max(\alpha,\beta)=\beta$ and $\max(\alpha,\gamma)=\gamma$. Therefore both sides of (\ref{E:minmax}) are $\sim$ to $\min(\beta,\gamma)$.

{\em Case 2.} Assume  $\min(\beta,\gamma)<\alpha$. Then $\max(\alpha,\min(\beta,\gamma))=\alpha$.  To decide  the right-hand side of (\ref{E:minmax}), suppose   $\beta\leq \gamma$ so we have the following  subcases.

($2a$) $\beta<\alpha\leq\gamma$: Then  $\max(\alpha,\beta)=\alpha$, $\max(\alpha,\gamma)=\gamma$, therefore, $\min(\max(\alpha,\beta),\max(\alpha,\gamma))=\alpha$, thus (\ref{E:minmax}) holds.

($2b$) $\beta\leq \gamma<\alpha$: Then $\max(\alpha,\beta)=\max(\alpha,\gamma)=\alpha$. So $$\min(\max(\alpha,\beta),\max(\alpha,\gamma))=\alpha,$$ thus (\ref{E:minmax}) holds again.

{\em Case 3.} Assume  $\min(\beta,\gamma)<\alpha$, so  $\max(\alpha,\min(\beta,\gamma))=\alpha$, but suppose now    $\gamma<\beta$. Then we have the subcases:

($3a$) $\gamma<\alpha\leq\beta$: Then $\max(\alpha,\beta)=\beta$ and $\max(\alpha,\gamma)=\alpha$. Thus $\min(\max(\alpha,\beta),\max(\alpha,\gamma))=\alpha$, that is,   (\ref{E:minmax}) holds.

($3b$) $\gamma< \beta\leq \alpha$: Then $\max(\alpha,\beta)=\alpha$ and $\max(\alpha,\gamma)=\alpha$. So $$\min(\max(\alpha,\beta),\max(\alpha,\gamma))=\alpha.$$ This completes the proof of the Proposition. \telos.

\vskip 0.2in

\begin{Cor} \label{C:eqschemes}
The schemes
\begin{equation} \label{E:scheme1}
\varphi\circ(\psi|\sigma)\leftrightarrow
(\varphi\circ\psi)|(\varphi\circ\sigma)
\end{equation}
(or its dual) and
\begin{equation} \label{E:scheme2}
\varphi\wedge\psi\rightarrow (\varphi|\psi\leftrightarrow\neg\varphi|\neg\psi)
\end{equation}
are equivalent and each one of them is a syntactic characterization of the regular $\neg$-decreasing orderings (and the corresponding choice functions).
\end{Cor}

{\em Proof.} The equivalence of (\ref{E:scheme1}) and (\ref{E:scheme2}), as {\em schemes}, follows from  Propositions \ref{P:notdistr},  \ref{P:distributive}, as well as from Lemma \ref{L:decrease} by which (\ref{E:scheme2}) characterizes the regular $\neg$-decreasing orderings. \telos

\vskip 0.2in

Interchanging $|$ and $\circ$ inside a sentence gives rise to a duality of sentences of $L_s$, that is, a mapping $\varphi\mapsto \varphi^d$ defined inductively as follows:

$\varphi^d=\varphi$, for classical $\varphi$,

$(\varphi\wedge\psi)^d=\varphi^d\wedge\psi^d$,

$(\neg\varphi)^d=\neg\varphi^d$,

$(\varphi|\psi)^d=\varphi^d\circ\psi^d$.

\vskip 0.1in

By the help of dual orderings $<^d$ and dual choice functions $f^d$, one can without much effort establish the following ``Duality Theorem'' which is the analogue of Boolean Duality:

\begin{Thm} \label{P:dualsame}
For every $\varphi\in Fml(L_s)$, $\models_{Reg^*}\varphi$ if and only if  $\models_{Reg^*}\varphi^d$.
\end{Thm}

\section{Axiomatization. Soundness and completeness results}

A {\em  Propositional Superposition Logic} (PLS for short) will consist as usual of two parts, a {\em syntactic} one, i.e., a formal system $K$, consisted of axiom-schemes and inference rules,  and a {\em semantical} one, consisted essentially of a set $X\subseteq {\cal F}$ of choice functions over $Sen(L)$.\footnote{More or less the same is true for every logical system, e.g.  PL. Although we  often  identify PL with the set of its logical axioms and the inference rule of Modus Ponens, tacitly we think of it  as a set of axiom-schemes $\textsf{Ax}({\rm PL})$ and  the inference rule $\textit{MP}$ on the ones hand, and the  natural Boolean semantics on the other. Specifically  $\textsf{Ax}({\rm PL})$ will consist of the following schemes:

(i) $\alpha\rightarrow (\beta\rightarrow\alpha)$

(ii) $(\alpha\rightarrow(\beta\rightarrow\gamma))
\rightarrow((\alpha\rightarrow\beta)\rightarrow(\alpha\rightarrow\gamma))$

(iii) $(\neg\alpha\rightarrow\neg\beta)\rightarrow ((\neg\alpha\rightarrow\beta)\rightarrow\alpha)$.
}
Let us start with the latter.

The semantical part $X$ induces the truth relation $\models_X$, that is, the class of structures $\langle M,f\rangle$, where $M:Sen(L)\rightarrow \{0,1\}$ and  $f\in X$, with respect to which the $X$-tautologies and $X$-logical consequence are defined.
For every  $X\subseteq {\cal F}$ let
$$Taut(X)=\{\varphi\in Sen(L_s):\models_X\varphi\}$$
be the set of tautologies of $L_s$ with respect to  $\models_X$. Let also $Taut$ be the set of classical tautologies. Then for any $X,Y\subseteq {\cal F}$,  $$X\subseteq Y \ \Rightarrow \ Taut\subseteq Taut(Y)\subseteq Taut(X).$$
The following simple fact will be used later but has also an obvious interest  in itself.

\begin{Lem} \label{L:computable}
For every $X\subseteq {\cal F}$, the set $Taut(X)$ is decidable (i.e., computable).
\end{Lem}

{\em Proof.} By the definition of $\models_s$,  $\varphi\in Taut(X)$ if and only if  $(\forall f\in X)(\overline{f}(\varphi)\in Taut)$ (where $Taut$ is the set of tautologies of PL), i.e.,
$$\varphi\in Taut(X) \ \Leftrightarrow \{\overline{f}(\varphi):f\in X\}\subset Taut.$$
Now given $\varphi$ and $f$, the collapse $\overline{f}(\varphi)$ results from $\varphi$ by inductively replacing each subformula $\psi_1|\psi_2$ of $\varphi$ with either $\overline{f}(\psi_1)$ or $\overline{f}(\psi_2)$. So clearly for every $\varphi$, the set of all possible collapses $\{\overline{f}(\varphi):f\in X\}$ is finite. Therefore, since $Taut$ is decidable, it is decidable whether $\{\overline{f}(\varphi):f\in X\}\subset Taut$.  \telos

\vskip 0.2in

In particular we are interested in the sets
\begin{equation} \label{E:inclusions}
Taut({\cal F})\subseteq  Taut(Reg)\subseteq Taut(Reg^*)\subseteq Taut(Dec),
\end{equation}
(as well as in $Taut(Asso)\subseteq Taut(Reg^*)$) corresponding to the truth relations considered above. It follows from \ref{L:computable} that these sets are decidable.  The question is whether each of  these sets of tautologies is  axiomatizable by a recursive set of axioms and inference rules. We shall see that the answer is yes.

Let us come to the formal system  $K$. Every $K$ consists of a set of axioms $\textsf{Ax}(K)$ and a set of inference rules $\textsf{IR}(K)$. Also the axioms and rules of  $K$ extend the axioms and rules of PL, i.e.,
$$\textsf{Ax}(K)=\textsf{Ax}({\rm PL})+\{S_i:i\leq n\}, \ \mbox{and} \ \textit{MP}\in \textsf{IR}(K),$$
where $S_i$ will be some schemes considered below expressing  basic properties of $|$. Given  $X\subseteq {\cal F}$ in order to axiomatize $Taut(X)$ by a formal system $K$, clearly it is necessary for the axioms of $K$ to be $X$-tautologies, i.e.,
$$\textsf{Ax}(K)\subseteq Taut(X).$$
For any such $X\subseteq {\cal F}$ and $K$, we  have a {\em logic}  that extends PL, called {\em   Propositional Superposition Logic w.r.t. to $X$ and $K$,} denoted $${\rm PLS}(X,K).$$

Given a formal system $K$ as above and $\Sigma\cup\{\varphi\}\subset Sen(L)$, a (Hilbert-style) {\em $K$-proof of $\varphi$ from $\Sigma$}   is defined just as a proof in  PL (mutatis mutandis), that is, as a sequence of sentences $\sigma_1,\ldots,\sigma_n$ such that $\sigma_n=\varphi$ and each $\sigma_i$ either belongs to $\Sigma$ or belongs to $\textsf{Ax}(K)$, or is derived from previous ones by the inference rules in   $\textsf{IR}(K)$.  We denote by
$$\Sigma\vdash_K \varphi$$
the fact that there is a $K$-proof of $\varphi$ from $\Sigma$. Especially for  classical  sentences, i.e.,  $\Sigma\cup\{\alpha\}\subseteq Sen(L)$, it is clear that
$$\Sigma\vdash_{PL} \alpha \  \Leftrightarrow \ \Sigma\vdash_K \alpha,$$
where $\vdash_{PL}$ denotes  provability in PL.
$\Sigma$ is said to be {\em $K$-consistent,}  if $\Sigma\not\vdash_K\bot$. Again for $\Sigma\subset Sen(L)$,
$$\Sigma \ \mbox{is $K$-consistent} \ \Leftrightarrow \ \Sigma \ \mbox{is consistent}.$$
Recall that a formal system $K$ (or its proof relation $\vdash_K$) satisfies the Deduction Theorem (DT) if for all $\Sigma$, $\varphi$ and $\psi$,
\begin{equation}
\Sigma\cup\{\varphi\}\vdash_K\psi \ \Rightarrow \ \Sigma\vdash_K\varphi\rightarrow\psi.
\end{equation}
It is well-known that if the only inference rule  of $K$ is $\textit{MP}$ (and perhaps also the  Generalization Rule), then  DT holds for $\vdash_K$. But in systems with additional inference rules DT often fails.  Below we shall consider formal  systems $K$  augmented with an additional  inference rule. So we shall need to examine the validity of DT later.

\begin{Def} \label{D:satifiable}
{\em A set   $\Sigma\subset Sen(L_s)$ is said to be}  $X$-satisfiable {\em  if for some truth assignment $M$ for $L$ and  some  $f\in X$, $\langle M,f\rangle\models_s\Sigma$.}
\end{Def}
As is well-known the Soundness and Completeness Theorems of a logic have two distinct formulations, which are not always equivalent, depending on the semantics and the validity of Deduction Theorem. For the logic ${\rm PLS}(X,K)$ these two forms, ST1 and ST2 for Soundness and CT1 and  CT2 for Completeness, are the following:

$$({\rm ST1}) \hspace{.5\columnwidth minus .5\columnwidth} \Sigma\vdash_K\varphi \ \Rightarrow \ \Sigma\models_X\varphi, \hspace{.5\columnwidth minus .5\columnwidth} \llap{}$$
$$({\rm ST2}) \hspace{.5\columnwidth minus .5\columnwidth}
\Sigma \ \mbox{is $X$-satisfiable} \ \Rightarrow \ \Sigma \ \mbox{is $K$-consistent}
\hspace{.5\columnwidth minus .5\columnwidth} \llap{}$$
$$({\rm CT1}) \hspace{.5\columnwidth minus .5\columnwidth} \Sigma\models_X\varphi \ \Rightarrow \ \Sigma\vdash_K\varphi,  \hspace{.5\columnwidth minus .5\columnwidth} \llap{}$$
$$({\rm CT2}) \hspace{.5\columnwidth minus .5\columnwidth}
\Sigma \ \mbox{is $K$-consistent} \ \Rightarrow \ \Sigma \ \mbox{is $X$-satisfiable}.
\hspace{.5\columnwidth minus .5\columnwidth} \llap{}$$

Concerning the relationship between ST1 and ST2 and between CT1 and CT2 for ${\rm PLS}(X,K)$ the following holds.

\begin{Fac} \label{F:eqsat}
(i) For every $X$
\begin{equation} \label{E:needeq}
\Sigma\not\models_X \varphi\ \Rightarrow \Sigma\cup\{\neg\varphi\} \
\mbox{is $X$-satisfiable}.
\end{equation}
As a consequence,  $({\rm ST1})\Leftrightarrow ({\rm ST2})$ holds for every ${\rm PLS}(X,K)$.

(ii) $({\rm CT1})\Rightarrow ({\rm CT2})$ holds for every ${\rm PLS}(X,K)$. If $\vdash_K$ satisfies DT, then the converse holds too, i.e., $({\rm CT1})\Leftrightarrow ({\rm CT2})$.
\end{Fac}

{\em Proof.} (i)  (\ref{E:needeq}) follows immediately from the definition of $\models_X$ and the fact that the truth is bivalent. Now $({\rm ST1})\Rightarrow ({\rm ST2})$ is straightforward. For the converse assume ${\rm ST2}$ and $\Sigma\not\models_X \varphi$. By (\ref{E:needeq}) $\Sigma\cup\{\neg\varphi\}$ is $X$-satisfiable. By ${\rm ST2}$, $\Sigma\cup\{\neg\varphi\}$ is $K$-consistent, therefore $\Sigma\not\vdash_K\varphi$.

(ii) $({\rm CT1})\Rightarrow ({\rm CT2})$ is also straightforward. For the converse assume CT2,  DT  and $\Sigma\not\vdash_K\varphi$. It is well-known that by DT the latter is equivalent to the $K$-consistency of $\Sigma\cup\{\neg\varphi\}$. By CT2, $\Sigma\cup\{\neg\varphi\}$ is $X$-satisfiable. Therefore $\Sigma\not\models_X\varphi$. \telos

\vskip 0.2in

In view of Fact \ref{F:eqsat} (i) we do not need to distinguish any more between ST1 and ST2, and can refer simply to ``sound'' logics.

However the distinction between CT1 and CT2 remains. This is also exemplified by considering the semantic analogue of DT. Given a class $X\subseteq {\cal F}$,  let us call   the implication:
\begin{equation} \label{E:sdt}
\Sigma\cup\{\varphi\}\models_X\psi \ \Rightarrow \ \Sigma\models_X\varphi\rightarrow \psi
\end{equation}
{\em Semantic Deduction Theorem for  $X$} (or, briefly,  SDT).  Here is a  relationship  between DT and SDT via CT1.

\begin{Fac} \label{F:dssd}
For every $X\subseteq {\cal F}$, SDT for $X$ is true. This implies that if the logic ${\rm PLS}(X,K)$ is sound and satisfies CT1, then $K$ satisfies DT.
\end{Fac}

{\em Proof.} That SDT holds for every $X\subseteq {\cal F}$ is an easily verified  consequence of the semantics $\models_X$. Now assume that  ${\rm PLS}(X,K)$ is sound (i.e., satisfies (equivalently) both ST1 and ST2),  satisfies CT1, and  $\Sigma\cup\{\varphi\}\vdash_K\psi$. By ST1  it follows that $\Sigma\cup\{\varphi\}\models_X\psi$. By SDT (\ref{E:sdt}) we have  $\Sigma\models_X\varphi\rightarrow\psi$. Then CT1 implies  $\Sigma\vdash_K\varphi\rightarrow\psi$, as required. \telos

\vskip 0.2in

Next we give a list of specific axiom-schemes (referred also to simply as {\em axioms}) about  $|$, certain nested groups of which are going to axiomatize  the truth relations $\models_{\cal F}$, $\models_{Reg}$, $\models_{Reg^*}$ and $\models_{Dec}$ considered  in the previous sections.

\vskip 0.1 in

($S_1$) $\varphi\wedge \psi\rightarrow \varphi|\psi$

($S_2$) $\varphi|\psi\rightarrow \varphi\vee\psi$

($S_3$) $\varphi|\psi\rightarrow \psi|\varphi$

($S_4$) $(\varphi|\psi)|\sigma\rightarrow \varphi|(\psi|\sigma)$

($S_5$) $\varphi\wedge\neg\psi\rightarrow (\varphi|\psi\leftrightarrow\neg\varphi|\neg\psi)$

\vskip 0.1in

We shall split the axiomatization of the four basic truth relations considered in the previous section in two parts. We shall consider first the basic truth relation $\models_{\cal F}$ relying on the entire class of functions ${\cal F}$, and then we shall consider  the rest stricter relations $\models_{Reg}$, $\models_{Reg^*}$ and $\models_{Dec}$. The reason is that the relation $\models_{\cal F}$ can be axiomatized by a formal system having  $MP$ as  the only inference rule, while the rest systems require formal systems augmented with a second rule. The latter requirement makes these systems considerably more complicated.

\subsection{Axiomatizing the truth relation $\models_{\cal F}$}
In this section we deal with the relation $\models_{\cal F}$ and show that it can be soundly and completely axiomatized by the first three axioms $S_1$-$S_3$ cited above and Modus Ponens ($MP$). We call this formal system $K_0$.   Namely
\begin{equation} \label{KO}
\textsf{Ax}(K_0)=\textsf{Ax}({\rm PL})+\{S_1,S_2,S_3\} \ \mbox{and} \ \ \textsf{IR}(K_0)=\{\textit{MP}\}.
\end{equation}
Observe that $S_1$ and $S_2$, combined  with the axioms of PL,  prove (in $K_0$)  $\varphi|\varphi\leftrightarrow \varphi$.

It is easy to see that the logic ${\rm PLS}({\cal F},K_0)$ is sound. Namely we have the following more general fact.

\begin{Thm} \label{T:sound}
Let $X\subseteq {\cal F}$. If  $K$ is a system such that   $\textsf{Ax}(K)\subset Taut(X)$ and $\textsf{IR}(K)=\{\textit{MP}\}$, then ${\rm PLS}(X,K)$ is sound.
\end{Thm}

{\em Proof.} Let $X$, $K$ be as stated and $\Sigma\vdash_K\varphi$.
Let $\varphi_1,\ldots,\varphi_n$, where $\varphi_n=\varphi$,  be a $K$-proof of $\varphi$. As usual we show that  $\Sigma\models_X\varphi_i$, for every $1\leq i\leq n$,  by induction on $i$. Given $i$, suppose the claim holds for all $j<i$, and let $\langle M,f\rangle\models_s\Sigma$, for some assignment $M$ and $f\in X$. We show that  $\langle M,f\rangle\models_s\varphi_i$. If $\varphi_i\in \Sigma$ this is obvious. If $\varphi_i\in \textsf{Ax}(K)$, then $\langle M,f\rangle\models_s\varphi_i$, because by assumption $\textsf{Ax}(K)\subset Taut(X)$ and $f\in X$. Otherwise, since $\textit{MP}$ is the only inference rule of $K$,  $\varphi_i$ follows by $\textit{MP}$ from sentences $\varphi_j$, $\varphi_k=(\varphi_j\rightarrow \varphi_i)$, for some $j,k<i$. By the induction assumption,  $\langle M,f\rangle\models_s\varphi_j$ and $\langle M,f\rangle\models_s\varphi_k$. Therefore $\langle M,f\rangle\models_s\varphi_i$. \telos

\vskip 0.2in

\begin{Cor} \label{C:soundko}
The logic ${\rm PLS}({\cal F},K_0)$ is sound.
\end{Cor}

{\em Proof.} By Theorem \ref{T:interpol} and Fact \ref{F:extend} (iv), $S_1$, $S_2$, $S_3$ are schemes that hold in $\langle M,f\rangle$ for all $f\in {\cal F}$, therefore $\textsf{Ax}(K_0)\subset Taut({\cal F})$. So the claim follows from \ref{T:sound}. \telos

\vskip 0.2in

{\bf Completeness of ${\rm PLS}({\cal F},K_0)$}
We come to the  completeness of the  logic ${\rm PLS}({\cal F},K_0)$. As usual, a  set $\Sigma\subseteq Sen(L_s)$ is said to be complete if for every $\varphi\in Sen(L_s)$, $\varphi\in \Sigma$ or $\neg\varphi\in \Sigma$. If $\Sigma$ is $K$-consistent and complete, then for every $\varphi\in Sen(L_s)$, $\varphi\in\Sigma\Leftrightarrow \neg\varphi\not\in\Sigma$. Moreover if $\Sigma\vdash_K\varphi$, then $\varphi\in \Sigma$.

Before coming to the logics introduced  in the previous subsection,  we shall give  a general satisfiability criterion. Fix a class  $X\subseteq{\cal F}$  of choice functions and a set of axioms  $K\subseteq Taut(X)$. Let  $\Sigma$ be a $K$-consistent and  complete set of sentences of  $L_s$ and let   $\Sigma_1=\Sigma\cap Sen(L)$ be the subset of $\Sigma$ that contains the classical sentences of $\Sigma$. Then clearly $\Sigma_1$ is a consistent and  complete  set of sentences of $L$. By the Completeness Theorem of PL, there exists a truth assignment $M$ for $L$ such that, for every $\alpha\in Sen(L)$
\begin{equation} \label{E:constant}
\alpha\in \Sigma_1 \ \Leftrightarrow \ M\models\alpha.
\end{equation}
Given $\Sigma$, $\Sigma_1$, $M$  satisfying (\ref{E:constant}), and  a set $X\subseteq{\cal F}$ of choice functions, the question is under what conditions $M$ can be paired with a function $f\in X$ such that $\langle M,f\rangle\models_s\Sigma$.   Below we give a  simple characterization of this fact which is the key characterization of  $X$-satisfiability.

\begin{Lem} \label{L:Henkin}
Let $X\subseteq {\cal F}$ and $K\subset Taut(X)$. Let also $\Sigma$ be a $K$-consistent and complete set of sentences of  $L_s$ and let  $\Sigma_1=\Sigma\cap Sen(L)$  and  $M$ satisfy  (\ref{E:constant}).  Then for every $f\in X$, $\langle  M,f\rangle\models\Sigma$ if and only if  for every $\varphi\in Sen(L_s)$,
\begin{equation} \label{E:kalisynthiki}
\varphi\in \Sigma \ \Rightarrow \overline{f}(\varphi)\in \Sigma.
\end{equation}
(Actually (\ref{E:kalisynthiki})  is equivalent to $$\varphi\in \Sigma \ \Leftrightarrow \overline{f}(\varphi)\in \Sigma,$$  but the other direction follows from (\ref{E:kalisynthiki}), the consistency and completeness of $\Sigma$ and the fact that $\overline{f}(\neg\varphi)=\neg\overline{f}(\varphi)$.)
\end{Lem}

{\em Proof.}  Pick an $f\in X$ and suppose  $\langle M,f\rangle\models_s\Sigma$. Then  by the completeness of $\Sigma$ and the definition of $\models_s$, for every $\varphi\in Sen(L_s)$,
$$\varphi\in \Sigma \ \Leftrightarrow \ \langle M,f\rangle\models_s\varphi \Leftrightarrow M\models \overline{f}(\varphi).$$
Now by (\ref{E:constant}), $M\models \overline{f}(\varphi)\Rightarrow \overline{f}(\varphi)\in \Sigma_1\subset\Sigma$. Therefore $\varphi\in \Sigma \ \Rightarrow \overline{f}(\varphi)\in \Sigma$.  Thus  (\ref{E:kalisynthiki}) holds.

Conversely, suppose (\ref{E:kalisynthiki}) is true.  To show that  $\langle M,f\rangle\models_s\Sigma$, pick some $\varphi\in \Sigma$.  By (\ref{E:kalisynthiki}) $\overline{f}(\varphi)\in \Sigma$. Then  $\overline{f}(\varphi)\in \Sigma_1$ since $\overline{f}(\varphi)$ is classical, so by  (\ref{E:constant}) $M\models \overline{f}(\varphi)$. This means that $\langle M,f\rangle\models_s\varphi$, as required.   \telos

\vskip 0.2in

We come next to the completeness of  ${\rm PLS}({\cal F},K_0)$. The essential step of the proof is the following Lemma.

\begin{Lem} \label{L:key}
Let $\Sigma$ be a $K_0$-consistent and  complete set of sentences of $L_s$. Then $\Sigma$ is ${\cal F}$-satisfiable.
\end{Lem}

{\em Proof.} Let $\Sigma$ be  $K_0$-consistent and  complete. Then for any $\varphi,\psi\in Sen(L_s)$, the possible subsets of $\Sigma$ whose elements are $\varphi|\psi$, $\varphi$, $\psi$ or  their negations  are the following:

\vskip 0.1in

(a1) $\{\varphi|\psi, \varphi, \psi\}\subset \Sigma$

(a2)  $\{\varphi|\psi, \varphi,\neg\psi\}\subset \Sigma$

(a3) $\{\varphi|\psi, \neg\varphi,\psi\}\subset \Sigma$

(a4) $\{\neg(\varphi|\psi),\neg\varphi,\neg\psi\}\subset \Sigma$

(a5)  $\{\neg(\varphi|\psi), \varphi,\neg\psi\}\subset \Sigma$

(a6) $\{\neg(\varphi|\psi), \neg\varphi,\psi\}\subset \Sigma$

\vskip 0.1in

\noindent The remaining cases,

(a7) $\{\varphi|\psi, \neg\varphi, \neg\psi\}\subset\Sigma$

(a8) $\{\neg(\varphi|\psi), \varphi, \psi\}\subset\Sigma$

\noindent are impossible because they contradict   $K_0$-consistency and completeness of $\Sigma$. Indeed,  in case (a7) we have  $\neg\varphi\wedge\neg\psi\in \Sigma$. Also  $\varphi|\psi\in\Sigma$, so by  $S_2$ and completeness, $\varphi\vee\psi\in \Sigma$, a contradiction. In  case (a8)  $\varphi\wedge\psi\in \Sigma$. Also  $\neg(\varphi|\psi)\in\Sigma$, so by $S_1$ and completeness  $\neg(\varphi\wedge\psi)\in\Sigma$, a contradiction.

Given a pair $\{\alpha,\beta\}$ we say that ``$\{\alpha,\beta\}$ satisfies   (ai)'' if for $\varphi=\alpha$ and $\psi=\beta$, the corresponding case (ai) above, for $1\leq i\leq 6$, holds. We define a choice function $g$ for $L$ as follows:
\begin{equation} \label{E:basicmap}
g(\alpha,\beta)=
\left\{\begin{array}{l}
               (i) \ \alpha, \ \mbox{if $\{\alpha,\beta\}$ satisfies  (a2) or (a6)}  \\
               (ii) \ \beta, \ \mbox{if  $\{\alpha,\beta\}$ satisfies  (a3) or (a5) } \\
               (iii) \ \mbox{any of the $\alpha$, $\beta$, if $\{\alpha,\beta\}$ satisfies (a1)
                or (a4).}
            \end{array} \right.
\end{equation}

\vskip 0.1in

{\em Claim.} $\overline{g}$ satisfies the implication (\ref{E:kalisynthiki}) of the previous Lemma.

\vskip 0.1in

{\em Proof of the Claim.} We prove (\ref{E:kalisynthiki})  by induction on the length of $\varphi$. For $\varphi=\alpha\in Sen(L)$, $\overline{g}(\alpha)=\alpha$, so (\ref{E:kalisynthiki}) holds trivially. Similarly the induction steps for $\wedge$ and $\neg$ follow immediately from the fact that $\overline{g}$ commutes with these connectives and the completeness of $\Sigma$. So the only nontrivial step of the induction is that for $\varphi|\psi$. It suffices to assume
\begin{equation} \label{E:as1}
\varphi\in \Sigma\ \Rightarrow \overline{g}(\varphi)\in\Sigma,
\end{equation}
\begin{equation} \label{E:as2}
\psi\in \Sigma\ \Rightarrow \overline{g}(\psi)\in\Sigma,
\end{equation}
and prove
\begin{equation} \label{E:symp}
\varphi|\psi\in \Sigma\ \Rightarrow \overline{g}(\varphi|\psi)\in\Sigma.
\end{equation}
Assume $\varphi|\psi\in \Sigma$. Then the only possible combinations of  $\varphi$, $\psi$ and their negations that can belong to $\Sigma$ are those of cases (a1), (a2) and (a3) above. To prove (\ref{E:symp}) it suffices to check that $\overline{g}(\varphi|\psi)\in\Sigma$ in each of these cases. Note that $\overline{g}(\varphi|\psi)=g(\overline{g}(\varphi),\overline{g}(\psi))=g(\alpha,\beta)$, where $\overline{g}(\varphi)=\alpha$ and $\overline{g}(\psi)=\beta$ are sentences of $L$, so (\ref{E:basicmap}) applies.

Case (a1):  Then  $\varphi\in \Sigma$ and  $\psi\in \Sigma$. By (\ref{E:as1}) and (\ref{E:as2}), $\overline{g}(\varphi)\in\Sigma$ and $\overline{g}(\varphi)\in\Sigma$. By definition (\ref{E:basicmap}), $\overline{g}(\varphi|\psi)=g(\overline{g}(\varphi),\overline{g}(\psi))$ can be either $\overline{g}(\varphi)$ or $\overline{g}(\psi)$. So in either case  $\overline{g}(\varphi|\psi)\in\Sigma$.

Case (a2):  Then $\varphi\in \Sigma$ and   $\neg\psi\in \Sigma$. By (\ref{E:as1}) and (\ref{E:as2}), $\overline{g}(\varphi)\in\Sigma$, $\overline{g}(\psi)\notin\Sigma$. Also  by  (\ref{E:basicmap}), $\overline{g}(\varphi|\psi)=g(\overline{g}(\varphi),\overline{g}(\psi))=\overline{g}(\varphi)$, thus $\overline{g}(\varphi|\psi)\in \Sigma$.

Case (a3):  Then  $\neg\varphi\in \Sigma$,  $\psi\in \Sigma$. By (\ref{E:as1}) and  (\ref{E:as2}), $\overline{g}(\varphi)\notin\Sigma$, $\overline{g}(\psi)\in\Sigma$. By (\ref{E:basicmap}), $\overline{g}(\varphi|\psi)=g(\overline{g}(\varphi),\overline{g}(\psi))=\overline{g}(\psi)$, thus $\overline{g}(\varphi|\psi)\in \Sigma$. This completes the proof of the Claim.

\vskip 0.1in

It follows that condition (\ref{E:kalisynthiki}) is true,  so by Lemma \ref{L:Henkin}, if $M\models \Sigma_1$ where $\Sigma_1=\Sigma\cap Sen(L)$, then  $\langle M,g\rangle\models\Sigma$, therefore $\Sigma$ is ${\cal F}$-satisfiable. \telos

\vskip 0.2in

Let us remark here that, since $\vdash_{K_0}$ satisfies the Deduction Theorem, by Fact \ref{F:eqsat} the two forms of completeness theorem CT1 and CT2 are equivalent for ${\rm PLS}({\cal F},K_0)$. So it is indifferent which one we are going to prove for the system  ${\rm PLS}({\cal F},K_0)$.

\begin{Thm} \label{T:MainC0}
{\rm (Completeness of ${\rm PLS}({\cal F},K_0$))} The logic ${\rm PLS}({\cal F},K_0)$ is complete.
That is, if $\Sigma$ is $K_0$-consistent, then $\Sigma$  is ${\cal F}$-satisfiable.
\end{Thm}

{\em Proof.} Let $\Sigma$ be $K_0$-consistent. Extend $\Sigma$ to a $K_0$-consistent and  complete  $\Sigma^*\supseteq\Sigma$.  By Lemma \ref{L:key},  $\Sigma^*$ is ${\cal F}$-satisfiable. Therefore so is $\Sigma$.  \telos

\begin{Cor} \label{C:decidable}
The set $\{\varphi: \ \vdash_{K_0}\varphi\}$ is decidable.
\end{Cor}

{\em Proof.} By the soundness and completeness of ${\rm PLS}({\cal F},K_0)$,
$\{\varphi: \ \vdash_{K_0}\varphi\}=Taut({\cal F})$.  But
$Taut({\cal F})$ is decidable by Lemma \ref{L:computable}.  \telos

\subsection{Axiomatizing the truth relations for the classes $Reg$, $Reg^*$ and $Dec$}
The next systems, $K_1$-$K_3$, are intended to capture in addition  the semantic property of regularity considered in section 2.3. We need to define $K_1$ so that  if $\vdash_{K_1}\varphi$ then $\varphi\in Taut(Reg)$, and vice-versa (if possible). Specifically, if  $\alpha\sim\alpha'$, we need  $K_1$ to prove, for every $\beta$, $\alpha|\beta\leftrightarrow \alpha'|\beta$, i.e.,   $$\vdash_{K_1}(\alpha|\beta\leftrightarrow \alpha'|\beta).$$
This cannot be captured by an axiom-scheme, since no scheme can express the relation $\sim$ of logical equivalence. It can be captured however by a new {\em inference rule.} Roughly we need a rule guaranteeing  that if $\varphi$, $\psi$ are logically equivalent, then $\varphi$ and $\psi$ can be interchanged {\em salva veritate} in expressions containing $|$, that is one entailing $\varphi|\sigma\leftrightarrow \psi|\sigma$, for every $\sigma$.\footnote{Of course substitution of logically equivalent sentences  salva veritate holds also in classical logic, that is, if $\alpha\sim\alpha'$ and $\alpha$ is a subformula of $\beta$, then $\beta[\alpha]\sim\beta[\alpha']$, where $\beta[\alpha']$ is the result of replacing $\alpha$ with $\alpha'$ within $\beta$. This however is a simple consequence of the compositional semantics of classical logic. In contrast the choice semantics of PLS is by no means  compositional.}
This is the following rule denoted $SV$ (for salva veritate):
$$(SV) \quad \quad \mbox{\em from} \ \ \varphi\leftrightarrow\psi \ \ \mbox{\em infer} \ \varphi|\sigma\leftrightarrow\psi|\sigma,$$
$$\quad \quad \quad \quad \quad \mbox{if} \ \varphi\leftrightarrow\psi \ \mbox{is provable in $K_0$}.$$
We see that  $SV$  is a ``conditional rule'', applied under constraints, much like the  Generalization Rule  of first-order logic (from $\varphi(x)$ infer $\forall x\varphi(x)$, if  $x$ is not free in the premises, and also the Necessitation Rule of modal logic (from $\varphi$ infer $\square\varphi$, if  $\vdash\varphi$).  It follows that  $SV$ operates according to the following:

\begin{Fac} \label{F:applySV}
Let $K$ be a formal system such that $SV\in \textsf{IR}(K)$.  If $\vdash_{K_0}(\varphi\leftrightarrow\psi)$, then   $\vdash_K(\varphi|\sigma\leftrightarrow\psi|\sigma)$ for every $\sigma$.
\end{Fac}

Note that since, according to Corollary \ref{C:decidable}, it is decidable, given   $\varphi$, whether $\vdash_{K_0}\varphi$, it is decidable, given a recursive set of sentences $\Sigma$,  whether a sequence of sentences $\varphi_1,\ldots,\varphi_n$ is a proof (in $K$) from $\Sigma$ or not.

\begin{Thm} \label{T:soundstrong}
Let $X\subseteq Reg$. If  $K$ is a system such that   $\textsf{Ax}(K)\subset Taut(X)$ and $\textsf{IR}(K)=\{\textit{MP},SV\}$, then ${\rm PLS}(X,K)$ is sound.
\end{Thm}

{\em Proof.}  Let  $X\subseteq Reg$, $\textsf{Ax}(K)\subset Taut(X)$ and $\textsf{IR}(K)=\{\textit{MP},SV\}$, and let $\Sigma\vdash_K\varphi$. Let $\varphi_1,\ldots,\varphi_n$, where $\varphi_n=\varphi$,  be a $K$-proof of $\varphi$. We show, by induction on $i$,  that for all $i=1,\ldots,n$,  $\Sigma\models_X\varphi_i$ . Let $\langle M,f\rangle\models_s\Sigma$, with $f\in X$.
The proof that $\langle M,f\rangle\models_s\varphi_i$ goes exactly as in the proof of Theorem \ref{T:sound}, except of the case where $\varphi_i$ follows from a sentence $\varphi_j$, for $j<i$, by $SV$. It means that $\varphi_i=(\sigma|\tau\leftrightarrow \rho|\tau)$ while $\varphi_j=(\sigma\leftrightarrow \rho)$, where $\vdash_{K_0}(\sigma\leftrightarrow \rho)$. Now $K_0$ is a system satisfying the conditions of \ref{T:sound} above for $X={\cal F}$,  so $\models_{\cal F}(\sigma\leftrightarrow \rho)$. It means that for every assignment $N$ and every $g\in {\cal F}$, $\langle N,g\rangle\models_s(\sigma\leftrightarrow \rho)$, i.e., $N\models\overline{g}(\sigma)\leftrightarrow \overline{g}(\rho)$, that is, $\overline{g}(\sigma)\leftrightarrow \overline{g}(\rho)$ is a classical tautology, or $\overline{g}(\sigma)\sim\overline{g}(\rho)$, for every $g\in {\cal F}$. In particular, $\overline{f}(\sigma)\sim\overline{f}(\rho)$. Now since $X\subseteq Reg$, $f\in X$ implies   $f$ is regular. Therefore   $\overline{f}(\sigma)\sim\overline{f}(\rho)$ implies that  $f(\overline{f}(\sigma),\overline{f}(\tau))\sim f(\overline{f}(\rho),\overline{f}(\tau))$, or $\overline{f}(\sigma|\tau)\sim \overline{f}(\rho|\tau)$, therefore $M\models \overline{f}(\sigma|\tau)\leftrightarrow \overline{f}(\rho|\tau)$, or $\langle M,f\rangle\models_s(\sigma|\tau\leftrightarrow \rho|\tau)$, i.e.,  $\langle M,f\rangle\models_s\varphi_i$, as required. This completes the proof. \telos

\vskip 0.2in

We define next the systems $K_1$-$K_3$ as follows:
\begin{equation} \label{E:k1}
\textsf{Ax}(K_1)=\textsf{Ax}(K_0)=\{S_1,S_2,S_3\},   \quad \ \quad \textsf{IR}(K_1)=\{\textit{MP},SV\},
\end{equation}
\begin{equation} \label{E:k2}
\textsf{Ax}(K_2)=\textsf{Ax}(K_1)+S_4,   \quad \ \quad \textsf{IR}(K_2)=\{\textit{MP},SV\},
\end{equation}
\begin{equation} \label{E:k3}
\textsf{Ax}(K_3)=\textsf{Ax}(K_2)+S_5,   \quad \   \quad \textsf{IR}(K_3)=\{\textit{MP},SV\}.
\end{equation}

\begin{Thm} \label{T:specSound}
{\rm (Soundness)} The logics   ${\rm PLS}(Reg,K_1)$, ${\rm PLS}(Reg^*,K_2)$ and  ${\rm PLS}(Dec,K_3)$ are sound.
\end{Thm}

{\em Proof.}  This follows essentially from the general soundness Theorem \ref{T:soundstrong}. We have $Dec\subseteq Reg^*\subseteq Reg$, so all these classes of choice functions satisfy the condition $X\subseteq Reg$ of \ref{T:soundstrong}. Also $\textsf{IR}(K_i)=\{\textit{MP},SV\}$, for $i=1,2,3$. By Corollary \ref{C:soundko}, $\textsf{Ax}(K_0)\subset Taut({\cal F})$, and $Taut({\cal F})\subseteq Taut(Reg)\subseteq Taut(Reg^*)\subseteq Taut(Dec)$. Since $\textsf{Ax}(K_1)=\textsf{Ax}(K_0)$, we have $\textsf{Ax}(K_1)\subseteq Taut(Reg)$, so it follows that ${\rm PLS}(Reg,K_1)$ is sound. Next  $\textsf{Ax}(K_2)=\textsf{Ax}(K_0)+S_4$, so to see that ${\rm PLS}(Reg^*,K_2)$ is sound, it suffices to see that $S_4\in Taut(Reg^*)$. But by  Theorem \ref{T:assoc}  $S_4\in Taut(Asso)\subset Taut(Reg^*)$. Therefore $\textsf{Ax}(K_2)\subseteq Taut(Reg^*)$, and we are done.
Finally $\textsf{Ax}(K_3)=\textsf{Ax}(K_2)+S_5$ and by  Theorem  \ref{T:syncar}, the  scheme $S_5$ characterizes  $\neg$-decreasingness, thus $S_5\in Taut(Dec)$. So $\textsf{Ax}(K_3)\subset Taut(Dec)$ and by \ref{T:soundstrong} ${\rm PLS}(Dec,K_3)$ is sound. \telos

\vskip 0.2in

The following Lemma will be essential for the completeness of the aforementioned logics, proved in the next section.

\begin{Lem} \label{L:compessential}
If $\Sigma\subset Sen(L_s)$ is closed with respect to $\vdash_{K_i}$, for some $i=1,2,3$,  and $\alpha,\alpha'$ are  sentences of $L$ such that $\alpha\sim\alpha'$, then for every $\beta$, $(\alpha|\beta\leftrightarrow \alpha'|\beta)\in\Sigma$.
\end{Lem}

{\em Proof.} Let $\alpha\sim\alpha'$. Then $\vdash_{PL}\alpha\leftrightarrow\alpha'$, hence also  $\vdash_{K_0}\alpha\leftrightarrow\alpha'$. By $SV\in \textsf{IR}(K_i)$ it follows that for every $\beta$, $\vdash_{K_i}\alpha|\beta\leftrightarrow\alpha'|\beta$. Therefore  $(\alpha|\beta\leftrightarrow \alpha'|\beta)\in\Sigma$ since $\Sigma$ is $\vdash_{K_i}$-closed. \telos

\begin{Quest} \label{Q:open}
Do the formal systems $K_1$-$K_3$ satisfy the Deduction Theorem ($DT$)?
\end{Quest}

\noindent We guess that the answer to this question is negative but we do not have a proof.   The standard way to prove $DT$ for $\vdash_{K_i}$ is to assume  $\Sigma\cup\{\varphi\}\vdash_{K_i}\psi$,  pick a proof $\psi_1,\ldots,\psi_n$ of $\psi$, with $\psi_n=\psi$, and show that $\Sigma\vdash_{K_i}\varphi\rightarrow \psi_i$, for every $i=1,\ldots,n$, by induction on $i$. The only crucial step of the induction is the one concerning  the rule $SV$, i.e., to show  that for any $\sigma,\sigma',\tau$, if $\Sigma\vdash_{K_i}\varphi\rightarrow (\sigma\leftrightarrow\sigma')$, and $\vdash_{K_0}(\sigma\leftrightarrow\sigma')$, then $\Sigma\vdash_{K_i}\varphi\rightarrow (\sigma|\tau\leftrightarrow\sigma'|\tau)$.
Now clearly
$$\vdash_{PL}(\sigma\leftrightarrow\sigma')\rightarrow (\varphi\rightarrow(\sigma\leftrightarrow\sigma')),$$
so also
$$\vdash_{K_0}(\sigma\leftrightarrow\sigma')\rightarrow (\varphi\rightarrow(\sigma\leftrightarrow\sigma')).$$
This combined with $\vdash_{K_0}(\sigma\leftrightarrow\sigma')$ and $\textit{MP}$ gives
$$\vdash_{K_0}\varphi\rightarrow(\sigma\leftrightarrow\sigma')),$$
and hence, by PL again,
$$\vdash_{K_0}(\varphi\rightarrow\sigma)\leftrightarrow
(\varphi\rightarrow\sigma').$$
By $SV$ it follows that
$$\Sigma\vdash_{K_i}((\varphi\rightarrow\sigma)|\tau)\leftrightarrow
((\varphi\rightarrow\sigma')|\tau).$$
However it is not clear if and how one can get  from the  latter the required derivation $\Sigma\vdash_{K_i}\varphi\rightarrow (\sigma|\tau\leftrightarrow\sigma'|\tau)$.

\vskip 0.1in

It follows from the preceding discussion that DT is open  for the formal systems  $K_i$, $i=1,2,3$. Now by Fact \ref{F:dssd}, if DT fails for $K_i$ then necessarily CT1 fails for the logics ${\rm PLS}(Reg,K_1)$, ${\rm PLS}(Reg^*,K_2)$ and ${\rm PLS}(Dec,K_3)$. This means that CT1 is also open  for the preceding  logics.
(In connection with  the  status of DT  note that, surprisingly enough,  the question about the validity of this theorem remains  essentially unsettled even  for  a logical theory  as old as modal logic, see \cite{HN12}.)

\vskip 0.2in

{\bf Completeness} We come to the completeness of the aforementioned logics based on the systems $K_1$-$K_3$. First in view of  the open status of DT for the systems  $K_1$-$K_3$  and Fact \ref{F:eqsat} (ii), we cannot identify  the two forms of completeness  CT1 and CT2 for these systems. We only know that (CT1)$\Rightarrow$(CT2). So  we  can hope to prove CT2 for $K_1$-$K_3$.

There is however another serious side-effect of the lack of DT. This is  that we don't know whether  every consistent set of sentences can be extended to a consistent and {\em complete} set. Clearly every consistent set $\Sigma$ can be extended (e.g. by Zorn's Lemma) to a {\em maximal} consistent set $\Sigma'\supseteq\Sigma$. But maximality of $\Sigma'$  cannot guarantee completeness without DT (while the converse is true). For, theoretically, $\Sigma'$ may be maximal consistent and yet there is a $\varphi$ such that $\varphi\notin \Sigma'$ and $\neg\varphi\notin\Sigma'$, in which case  $\Sigma'\cup\{\varphi\}$ and $\Sigma'\cup\{\neg\varphi\}$ are both inconsistent. That  looks strange but we don't see how it could be proved false without DT. This property of extendibility of a consistent set to a consistent and complete one, for a formal system $K$, plays a crucial  role in the proof of completeness of $K$ (with respect to a given semantics), so we isolate it as property of $K$ denoted $cext(K)$. Namely we set
$$(cext(K))
\hspace{.5\columnwidth minus .5\columnwidth} \mbox{\em Every $K$-consistent set of sentences can be extended to }\hspace{.5\columnwidth minus .5\columnwidth} \llap{}$$
\hspace{0.9in} $\mbox{\em a $K$-consistent and complete set}$.

\vskip 0.1in

In view of the unknown truth-value of $cext(K_i)$, for $i=1,2,3$, we shall prove only {\em conditional} versions of ${\rm CT2}$-completeness for these systems. Actually it is shown that   ${\rm CT2}$-completeness is {\em equivalent} to $cext(K_i)$.

\begin{Thm} \label{T:MainC1}
{\rm (Conditional CT2-completeness  for ${\rm PLS}(Reg,K_1$))} The logic ${\rm PLS}(Reg,K_1$) is {\rm CT2}-complete if and only if $cext(K_1)$ is true.
\end{Thm}

{\em Proof.} One direction is easy. Assume $cext(K_1)$ is false. Then there is a maximal $K_1$-consistent set of sentences  $\Sigma$ non-extendible to a $K_1$-consistent and complete one. It means that there is a sentence $\varphi$ such that both $\Sigma\cup\{\varphi\}$ and $\Sigma\cup\{\neg\varphi\}$ are $K_1$-inconsistent. But then $\Sigma$ is not $Reg$-satisfiable. For if there are $M$ and $f\in Reg$ such that $\langle M,f\rangle\models_s\Sigma$, then $\langle M,f\rangle$ satisfies also either $\varphi$ or $\neg\varphi$. Thus either $\Sigma\cup\{\varphi\}$ or  $\Sigma\cup\{\neg\varphi\}$ is $Reg$-satisfiable. But this is a contradiction since both $\Sigma\cup\{\varphi\}$ and $\Sigma\cup\{\neg\varphi\}$ are inconsistent and by Theorem \ref{T:sound} ${\rm PLS}(Reg,K_1)$ is sound. Therefore $\Sigma$ is consistent and not $Reg$-satisfiable, so ${\rm PLS}(Reg,K_1)$ is not CT2-complete.

We come to the main direction of the equivalence assuming $cext(K_1)$ is true. Then given a $K_1$-consistent set $\Sigma$, we may assume without loss of generality that it is also complete.  We have to find $M$ and $g\in Reg$ such that $\langle M,g\rangle\models_s\Sigma$. It  turns out that the main argument of Lemma \ref{L:key}, concerning the definition of the choice  function $g$, works  also, with the necessary adjustments, for the other logics defined in the previous section.
Namely it suffices  to find  a  choice function $g\in Reg$  such that  $\langle M,g\rangle\models\Sigma$, where $M$ is a model of $\Sigma_1=\Sigma\cap Sen(L)$. The definition of  $g$ follows exactly the pattern of definition of $g$ in the proof of Lemma \ref{L:key}, except that we need now to take care  so that $g$ be regular. Recall that $g$ is regular if for all $\alpha$, $\alpha'$, $\beta$, $$\alpha'\sim\alpha \ \Rightarrow \ g(\alpha',\beta)\sim g(\alpha,\beta).$$ In (\ref{E:basicmap})  $g$ is defined by three clauses: (i) (a2) or (a6), (ii) (a3) or (a5), (iii) (a1) or (a4).

\vskip 0.1in

{\em Claim. } The regularity constraint is satisfied whenever $g$ is defined by clauses (i) and (ii) above.

\vskip 0.1in

{\em Proof of Claim.} Pick $\alpha$, $\alpha'$, $\beta$ such that $\alpha\sim\alpha'$. We prove the Claim for the case that  $g(\alpha,\beta)$ is defined according to  clause (i)-(a2). All other cases are verified  similarly.  That $g(\alpha,\beta)$  is defined by case (i)-(a2) of (\ref{E:basicmap})  means that $\alpha|\beta\in \Sigma$, $\alpha\in \Sigma$, $\neg\beta\in \Sigma$ and $g(\alpha,\beta)=\alpha$.  It suffices to see that necessarily $g(\alpha',\beta)=\alpha'\sim g(\alpha,\beta)$.

Since $\Sigma$ is complete, it is closed with respect to $\vdash_{K_1}$, so by Lemma \ref{L:compessential}, $\alpha\sim\alpha'$ implies that $(\alpha|\beta\leftrightarrow \alpha'|\beta)\in\Sigma$.  Also by assumption, $\alpha|\beta\in \Sigma$, hence $\alpha'|\beta\in \Sigma$. Moreover $\alpha'\in \Sigma$, since $\alpha\in \Sigma$, and $\neg\beta\in\Sigma$. Therefore case (i)-(a2) occurs too for $\alpha'|\beta$, $\alpha'$ and $\beta$. So, by (\ref{E:basicmap}),  $g(\alpha',\beta)=\alpha'$, therefore $g(\alpha',\beta)\sim g(\alpha,\beta)$. This proves the Claim.

\vskip 0.1in

It follows from the Claim that if we define  $g$ according to  (\ref{E:basicmap}), regularity  is  guaranteed unless   $g(\alpha,\beta)$ is given by  clause (iii), that is, unless  (a1) or (a4) is the case. In such a case either both $\alpha$, $\beta$ belong to $\Sigma$, or both $\neg\alpha$, $\neg\beta$ belong to $\Sigma$, and   (\ref{E:basicmap}) allows $g(\alpha,\beta)$ to be {\em any}  of the elements $\alpha$, $\beta$. So at this point  we must intervene by a new condition that will  guarantee regularity. This is done as follows.

Pick, as in the proof of Proposition \ref{P:ac}, from each $\sim$-equivalence class $[\alpha]$, a representative $\xi_\alpha\in[\alpha]$. Recall that, by completeness, the set $\Sigma_1=\Sigma\cap Sen(L)$ as well as its complement $\Sigma_2=Sen(L)-\Sigma_1$ are saturated with respect to $\sim$, that is, for every $\alpha$, either $[\alpha]\subset \Sigma_1$ or $[\alpha]\subset \Sigma_2$. Let $D_1=\{\xi_\alpha:\alpha\in \Sigma_1\}$, $D_2=\{\xi_\alpha:\alpha\in \Sigma_2\}$. Let $[D_i]^2$ be the set of pairs of elements  of $D_i$, for $i=1,2$, and pick  an arbitrary choice function $g_0:[D_1]^2\cup [D_2]^2\rightarrow D_1\cup D_2$.
Then it suffices to define $g$ by slightly revising definition  (\ref{E:basicmap}) as follows:

\begin{equation} \label{E:basicmap1}
g(\alpha,\beta)=
\left\{\begin{array}{l}
               (i) \ \alpha, \ \mbox{if $\{\alpha,\beta\}$,  satisfies  (a2) or (a6)}  \\
               (ii) \ \beta, \ \mbox{if  $\{\alpha,\beta\}$ satisfies  (a3) or (a5) } \\
               (iii) \ \sim g_0(\xi_{\alpha},\xi_{\beta}), \  \mbox{if $\{\alpha,\beta\}$ satisfies (a1) or (a4).}
            \end{array} \right.
\end{equation}
(The third clause is  just a shorthand for: $g(\alpha,\beta)=\alpha$ if $g_0(\xi_{\alpha},\xi_{\beta})=\xi_{\alpha}$, and $g(\alpha,\beta)=\beta$ if $g_0(\xi_{\alpha},\xi_{\beta})=
\xi_{\beta}$. See the similar formulation in the proof of \ref{P:ac}.)
In view of the Claim and the specific definition of $g$ by (\ref{E:basicmap1}), it follows immediately that  if $\alpha\sim \alpha'$ then for every  $\beta$, $g(\alpha,\beta)\sim g(\alpha',\beta)$.  So $g$ is regular. Further,  exactly  as in Lemma \ref{L:key} it follows that $\langle M,g\rangle\models_s\Sigma$. This completes the proof.  \telos

\vskip 0.2in

Next we come to the logic ${\rm PLS}(Reg^*,K_2)$. The difference of  $K_2$-consistency from  $K_1$-consistency is that, as a result of axiom $S_4$, if $\Sigma$ is $K_2$-consistent and $\varphi|(\psi|\sigma)\in \Sigma$, then $(\varphi|\psi)|\sigma\in \Sigma$, or more simply $\varphi|\psi|\sigma\in \Sigma$. Let us outline this difference by an example.

\begin{Ex} \label{Ex:counter}
 Let
$$\Sigma=\{\alpha, \neg\beta, \neg\gamma,\alpha|\beta,
\neg(\alpha|\gamma), \alpha|(\beta|\gamma)\},$$
where  $\alpha$, $\beta$, $\gamma$ are pairwise inequivalent and $\alpha\wedge\neg\beta\wedge\neg\gamma$ is satisfiable. Then $\Sigma$ is $Reg$-satisfiable, hence $K_1$-consistent, but is not $K_2$-consistent.
In particular $\Sigma$ is not $Asso$-satisfiable.
\end{Ex}

{\em Proof.} By hypothesis there is a truth assignment $M$ such that $M\models \alpha\wedge\neg\beta\wedge\neg\gamma$. Pick a (partial) choice function for $L$ such that $f(\alpha,\beta)=\alpha$,  $f(\alpha,\gamma)=\gamma$ and $f(\beta,\gamma)=\beta$. Since $\alpha$, $\beta$, $\gamma$ are pairwise inequivalent, it is easy to see that $f$ extends to a regular choice function for the entire $L$. Then $\overline{f}(\alpha|\beta)=\alpha$,   $\overline{f}(\alpha|\gamma)=\gamma$ and $\overline{f}(\beta|\gamma)=\beta$.  So  $\overline{f}(\neg(\alpha|\gamma))=\neg\gamma$. It follows that $\langle M,f\rangle\models_s \{\alpha|\beta,\neg(\alpha|\gamma)\}$. Moreover $\overline{f}(\alpha|(\beta|\gamma))=f(\alpha,f(\beta,\gamma))=
f(\alpha,\beta)=\alpha$, which means that   $\langle M,f\rangle\models_s \alpha|(\beta|\gamma)$ too. Thus $\langle M,f\rangle\models_s\Sigma$, so $\Sigma$ is $Reg$-satisfiable.

Now in view of axiom $S_4$ of $K_2$, since $\alpha|(\beta|\gamma)\in\Sigma$ it follows that $\Sigma\vdash_{K_2}\alpha|\beta|\gamma$. By $S_2$ and $S_3$, the latter implies   $\Sigma\vdash_{K_2}(\alpha|\gamma)\vee\beta$. On the other hand $\neg(\alpha|\gamma)\in\Sigma$ and $\neg\beta\in\Sigma$, so  $\Sigma\vdash_{K_2}\neg(\alpha|\gamma)\wedge \neg\beta$, or $\Sigma\vdash_{K_2}\neg((\alpha|\gamma)\vee\beta)$. Thus $\Sigma\vdash_{K_2}\bot$, so  it is $K_2$-inconsistent.

Finally  assume that $\Sigma$ is satisfied in $\langle N,f\rangle$, for some assignment $N$ and some associative $f$. Let $f=\min_<=\min$ for some total ordering $<$ of $Sen(L)$. Now $\langle N,f\rangle\models_s\{\alpha, \neg\beta,\alpha|\beta\}$ implies $\min(\alpha,\beta)=\alpha$, i.e., $\alpha<\beta$, while $\langle N,f\rangle\models_s\{\alpha, \neg\gamma,\neg(\alpha|\gamma)\}$ implies $\min(\alpha,\gamma)=\gamma$, so  $\gamma<\alpha$. Therefore $\gamma<\alpha<\beta$. On the other hand, $\langle N,f\rangle\models_s \alpha|(\beta|\gamma)$ implies $N\models\min(\alpha,\min(\beta,\gamma))=\min(\alpha,\beta,\gamma)$, therefore  $\min(\alpha,\beta,\gamma)=\alpha$ since $N\models\neg\beta\wedge\neg\gamma$. Thus $\alpha<\gamma$, a contradiction. \telos

\begin{Thm} \label{T:MainC2}
{\rm (Conditional CT2-completeness  for ${\rm PLS}(Reg^*,K_2$))} The logic ${\rm PLS}(Reg^*,K_2$) is {\rm CT2}-complete if and only if $cext(K_2)$ is true.
\end{Thm}

{\em Proof.} One direction of the equivalence is proved exactly as the  corresponding direction  of Theorem \ref{T:MainC1}. So let us come to the other direction assuming $cext(K_2)$ is true. Let $\Sigma$ be a $K_2$-consistent set, so we may  assume again that $\Sigma$ is also complete. We must construct a  regular and associative choice function $g$  such that $\langle M,g\rangle\models\Sigma$, where $M\models \Sigma_1$. As already remarked,  $\alpha|(\beta|\gamma)\in \Sigma$ implies $(\alpha|\beta)|\gamma\in \Sigma$. We shall define   $\overline{g}$  basically as in  definition (\ref{E:basicmap}) of Lemma \ref{L:key}, except that now we want $g$ to induce  a regular total ordering of $Sen(L)$. So let $h$ be a partial choice function for $L$ such that
$$dom(h)=\{\{\alpha,\beta\}: \{\alpha,\beta\}
\ \mbox{satisfies some of the cases (a2), (a3), (a5) and}$$
\hspace{0.8in} (a6) of Lemma \ref{L:key}\}, \\
and
\begin{equation} \label{E:thiish}
h(\alpha,\beta)=
\left\{\begin{array}{l}
               (i) \ \alpha, \ \mbox{if $\{\alpha,\beta\}$ satisfies  (a2) or (a6),}  \\
               (ii) \ \beta, \ \mbox{if  $\{\alpha,\beta\}$ satisfies  (a3) or (a5). } \\
            \end{array} \right.
\end{equation}

\vskip 0.1in

{\em Claim 1.} For any $\alpha$, $\beta$, $\gamma$, whenever two of  the $h(\alpha,h(\beta,\gamma))$, $h(\beta, h(\alpha,\gamma))$, $h(\gamma, h(\alpha,\beta))$ are defined,  they are equal.

\vskip 0.1in

{\em Proof of Claim 1.} Pick some $\alpha$, $\beta$, $\gamma$. Then at least two of them belong  either to $\Sigma$ or to its complement. Without loss of generality assume  that $\alpha\in \Sigma$, $\beta\notin \Sigma$, $\gamma\notin \Sigma$. Then in view of $K_2$-consistency and completeness of $\Sigma$,  $$A=\{\alpha,\neg\beta,\neg\gamma, \neg(\beta|\gamma)\}\subset\Sigma.$$
Also by $K_2$-consistency and completeness we can identify $\alpha|(\beta|\gamma)$ and  $(\alpha|\beta)|\gamma$, with respect to their containment to $\Sigma$, and there are two options: either $\alpha|\beta|\gamma\in \Sigma$ or $\neg(\alpha|\beta|\gamma)\in \Sigma$.
We consider now the  combinations of the sentences $\alpha|\beta|\gamma$, $\alpha|\beta$, $\alpha|\gamma$ and their negations  that can belong to $\Sigma$ together with the elements of $A$. We write these combinations  in the form of sets $B_i$. It is easy to see that the only sets $B_i$ of this kind such that  $A\cup B_i\subset \Sigma$, are the following:

$B_1=\{\alpha|\beta|\gamma, \alpha|\beta,\alpha|\gamma\}$

$B_2=\{\neg(\alpha|\beta|\gamma), \alpha|\beta,\alpha|\gamma\}$

$B_3=\{\neg(\alpha|\beta|\gamma), \alpha|\beta,\neg(\alpha|\gamma)\}$

$B_4=\{\neg(\alpha|\beta|\gamma), \neg(\alpha|\beta),\alpha|\gamma\}$

$B_5=\{\neg(\alpha|\beta|\gamma), \neg(\alpha|\beta),\neg(\alpha|\gamma)\}$

\vskip 0.1in

\noindent The remaining sets:

$B_1'=\{\alpha|\beta|\gamma, \alpha|\beta,\neg(\alpha|\gamma)\}$

$B_2'=\{\alpha|\beta|\gamma, \neg(\alpha|\beta),\alpha|\gamma\}$

$B_3'=\{\alpha|\beta|\gamma, \neg(\alpha|\beta),\neg(\alpha|\gamma)\}$

\noindent cannot be included in $\Sigma$ jointly with $A$. [For instance assume  $$A\cup\{\alpha|\beta|\gamma,\alpha|\beta,\neg(\alpha|\gamma)\}
=\{\alpha,\neg\beta,\neg\gamma, \neg(\beta|\gamma),\alpha|\beta|\gamma,\alpha|\beta,\neg(\alpha|\gamma)\}
\subset\Sigma.$$
Then $\alpha|\beta|\gamma$ is written $\beta|\alpha|\gamma$ so by $S_2$, it implies $\beta\vee (\alpha|\gamma)\in \Sigma$. By completeness, either $\beta\in \Sigma$ or $\alpha|\gamma\in \Sigma$. But already  $\neg\beta$ and $\neg(\alpha|\gamma)$ are in $\Sigma$, a contradiction.]

Since $\beta,\gamma\notin\Sigma$, $h(\beta,\gamma)$, and hence $h(\alpha,h(\beta,\gamma))$ are not defined by (i) or (ii) of (\ref{E:basicmap}). So it suffices to verify that $h(\beta, h(\alpha,\gamma))=h(\gamma, h(\alpha,\beta))$ in each of the cases $A\cup B_i\subset \Sigma$, for $1\leq i\leq 5$.

1) $A\cup B_1\subset\Sigma$: We have $h(\alpha,\beta)=\alpha$, $h(\alpha,\gamma)=\alpha$. Then
$$h(\beta, h(\alpha,\gamma))=h(\beta,\alpha)=\alpha=h(\alpha,\gamma)=h(\gamma, h(\alpha,\beta)),$$
so the Claim holds.

2) $A\cup B_2\subset\Sigma$:  Same as before.

3) $A\cup B_3\subset\Sigma$: We have $h(\alpha,\beta)=\alpha$, $h(\alpha,\gamma)=\gamma$. Thus $h(\beta,h(\alpha,\gamma)=h(\beta,\gamma)$, so $h(\beta,h(\alpha,\gamma)$ is also undefined. We se that only $h(\gamma,h(\alpha,\beta))=\gamma$ is defined, so  the Claim holds vacuously.

4) $A\cup B_4\subset\Sigma$: We have $h(\alpha,\beta)=\beta$, $h(\alpha,\gamma)=\alpha$. Thus  $h(\gamma,h(\alpha,\beta))=g(\gamma,\beta)$ is undefined, and the  Claim holds vacuously as before.

5) $A\cup B_5\subset\Sigma$: We have $h(\alpha,\beta)=\beta$, $h(\alpha,\gamma)=\gamma$. Thus $h(\beta,h(\alpha,\gamma))=g(\beta,\gamma)$ is undefined, and we are done again. This  completes  the proof of  Claim 1.

\vskip 0.1in

{\em Claim 2.} Let   $$S=\{\langle \alpha,\beta\rangle: \{\alpha,\beta\}\in dom(h) \wedge  h(\alpha,\beta)=\alpha\},$$ and let $<_1$ be the transitive closure of $S$. Then $<_1$ is a regular partial ordering on $Sen(L)$.

\vskip 0.1in

{\em Proof of Claim 2.} By Claim 1, $h$ is an associative partial choice function, so as in the proof of Theorem \ref{T:LO} we can see that the transitive closure $<_1$ of $S$  is a partial ordering. Also that $<_1$ is a regular ordering  follows from the Claim of Theorem \ref{T:MainC1}.  This  completes  the proof of  Claim 2.

\vskip 0.1in

Now clearly the partial ordering $<_1$ of Claim 2 extends to a regular total ordering $<$ of $Sen(L)$. Then it suffices to define $g$ by setting $g=\min_<$. Since  for every $\alpha,\beta\in Sen(L)$, if the pair $\langle \alpha,\beta\rangle$ satisfies some of the cases (a2), (a3), (a5), (a6), $\alpha<\beta$ if and only if  $h(\alpha,\beta)=\alpha$, clearly $g$ extends $h$. Moreover

\begin{equation} \label{E:basicmap2}
g(\alpha,\beta)=
\left\{\begin{array}{l}
               (i) \ \alpha, \ \mbox{if $\langle\alpha,\beta\rangle$ satisfies  (a2) or (a6),}  \\
               (ii) \ \beta, \ \mbox{if  $\langle\alpha,\beta\rangle$ satisfies  (a3) or (a5), } \\
               (iii) \ \min_<(\alpha,\beta), \  \mbox{if $\langle\alpha,\beta\rangle$ satisfies (a1) or (a4).}
            \end{array} \right.
\end{equation}
Thus  it follows as in Lemma \ref{L:key} that $\langle M,g\rangle\models_s\Sigma$, that is, $\Sigma$ is $Reg^*$-satisfiable. This completes the proof of the theorem. \telos

\vskip 0.2in

Finally we come to the conditional completeness  of ${\rm PLS}(Dec,K_3)$.

\begin{Thm} \label{T:MainC3}
{\rm (Conditional CT2-completeness  for ${\rm PLS}(Dec,K_3$))} The logic ${\rm PLS}(Dec,K_3$) is {\rm CT2}-complete if and only if $cext(K_3)$ is true.
\end{Thm}

{\em Proof.} Again one direction of the equivalence is shown  exactly as the  corresponding direction  of Theorem \ref{T:MainC1}. We come to the other direction assuming $cext(K_3)$ is true. Fix a  $K_3$-consistent set $\Sigma$. By $cext(K_3)$ we may assume that  $\Sigma$ is also complete. Let $M\models\Sigma_1$, where $\Sigma_1=\Sigma\cap Sen(L)$. We show that there exists a choice function $g$ such that $g=\min_<$, where $<$ is a  $\neg$-decreasing regular total ordering of $Sen(L)$, and $\langle M,g\rangle\models \Sigma$. $g$ is essentially defined as in the previous theorem plus an extra adjustment that guarantees $\neg$-decreasingness. Namely, let $h$ be the function defined exactly as in the proof of Theorem \ref{T:MainC2}.

\vskip 0.1in

{\em Claim.} $h$ is $\neg$-decreasing, i.e., whenever $h(\alpha,\beta)$ and $h(\neg\alpha,\neg\beta)$ are defined, then
\begin{equation} \label{E:showmono}
h(\alpha,\beta)=\alpha\ \Leftrightarrow h(\neg\alpha,\neg\beta)=\neg\beta.
\end{equation}

\vskip 0.1in

{\em Proof of  Claim.} We must check that whenever $\{\alpha,\beta\}$ and $\{\neg\alpha,\neg\beta\}$ satisfy some of the cases  (a2), (a3), (a5) and (a6), then (\ref{E:showmono}) holds true. Thus we must examine the   combinations of $\alpha$, $\beta$, $\alpha|\beta$, $\neg\alpha|\neg\beta$ and their negations that can belong to $\Sigma$. There is a total of 16 possible combinations of these sentences. Of them the combinations

\vskip 0.1in

$U_1=\{\alpha|\beta, \neg(\neg\alpha|\neg\beta), \alpha,\beta\}$

$U_2=\{\neg(\alpha|\beta), \neg\alpha|\neg\beta,\neg\alpha,\neg\beta\}$

\vskip 0.1in

\noindent do not allow definition of $h$ since in these cases either both $\alpha$, $\beta$ or both $\neg\alpha$, $\neg\beta$ belong to $\Sigma$.
Next we have 10   combinations that contradict $K_3$-consistency and completeness of $\Sigma$. These are:

\vskip 0.1in

$F_1=\{\alpha|\beta, \neg\alpha|\neg\beta, \alpha,\beta\}$

$F_2=\{\alpha|\beta, \neg\alpha|\neg\beta, \neg\alpha,\neg\beta\}$

$F_3=\{\alpha|\beta, \neg(\neg\alpha|\neg\beta), \neg\alpha,\neg\beta\}$

$F_4=\{\alpha|\beta, \neg(\neg\alpha|\neg\beta), \neg\alpha,\beta\}$

$F_5=\{\alpha|\beta, \neg(\neg\alpha|\neg\beta), \alpha,\neg\beta\}$

$F_6=\{\neg(\alpha|\beta), \neg\alpha|\neg\beta, \alpha,\beta\}$

$F_7=\{\neg(\alpha|\beta), \neg\alpha|\neg\beta, \neg\alpha,\beta\}$

$F_8=\{\neg(\alpha|\beta), \neg\alpha|\neg\beta, \alpha,\neg\beta\}$

$F_9=\{\neg(\alpha|\beta), \neg(\neg\alpha|\neg\beta), \alpha,\beta\}$

$F_{10}=\{\neg(\alpha|\beta), \neg(\neg\alpha|\neg\beta), \neg\alpha,\neg\beta\}$.

\vskip 0.1in

\noindent Notice  that of the preceding sets, $F_4$, $F_5$, $F_7$ and  $F_8$ yield a contradiction because of the axiom $S_5$. For instance consider $F_4=\{\alpha|\beta, \neg(\neg\alpha|\neg\beta), \neg\alpha,\beta\}$. It  contains $\neg\alpha,\beta$, thus it proves $\neg\alpha\wedge\beta$. By  $S_5$,  $F_4$ proves $(\alpha|\beta\leftrightarrow \neg\alpha|\neg\beta)$. $F_4$ also contains $\alpha|\beta$, thus it proves $\neg\alpha|\neg\beta$. But it besides  contains $\neg(\neg\alpha|\neg\beta)$, so $F_4\vdash_{K_3}\bot$.
Thus the only  combinations that can be contained in $\Sigma$ are the following:

\vskip 0.1in

$C_1=\{\alpha|\beta, \neg\alpha|\neg\beta, \alpha,\neg\beta\}\subset \Sigma$

$C_2=\{\alpha|\beta, \neg\alpha|\neg\beta, \neg\alpha,\beta\}\subset \Sigma$

$C_3=\{\neg(\alpha|\beta), \neg(\neg\alpha|\neg\beta), \alpha,\neg\beta\}\subset\Sigma$

$C_4=\{\neg(\alpha|\beta), \neg(\neg\alpha|\neg\beta), \neg\alpha,\beta\}\subset\Sigma$.

\vskip 0.1in

\noindent It is easy to verify that in each of the  cases $C_i\subset \Sigma$, for $1\leq i\leq 4$, (\ref{E:showmono}) is true in view of the definition (\ref{E:thiish}) of $h$. For example in case $C_4\subset \Sigma$, necessarily $h(\alpha,\beta)=\alpha$, while $h(\neg\alpha,\neg\beta)=\neg\beta$. This completes the proof of  Claim 1.

\vskip 0.1in

As in the proof of \ref{T:MainC2}, let $$S=\{\langle \alpha,\beta\rangle:\{\alpha,\beta\}\in dom(h) \wedge h(\alpha,\beta)=\alpha\},$$ and let $<_1$ be the transitive closure of $S$. As shown in \ref{T:MainC2}, $<_1$ is a regular partial ordering. Moreover here, in view of the Claim, $<_1$ is $\neg$-decreasing. So by a standard application of Zorn's Lemma, $<_1$ extends to a regular $\neg$-decreasing total ordering $<$ of $Sen(L)$. If we set $g=\min_<$, then $g$ satisfies  (\ref{E:basicmap2}) of the previous theorem and thus $\langle M,g\rangle\models_s\Sigma$. Therefore  $\Sigma$ is $Dec$-satisfiable.   \telos

\vskip 0.2in

The following is open.

\begin{Quest} \label{Q:CT1}
If  $cext(K_i)$ are true for $i=1,2,3$, do the logics ${\rm PLS}(Reg,K_1)$, ${\rm PLS}(Reg^*,K_2)$ and ${\rm PLS}(Dec,K_3)$ satisfy the form {\rm CT1} of Completeness Theorem?
\end{Quest}

\subsection{Some closing remarks on axiomatization}

Before closing this section on axiomatization of superposition logics, let us notice that all axioms $S_1$-$S_5$ introduced above are true also for the connectives $\wedge$ and $\vee$. That is, none of the $S_i$ can be used to discriminate $|$ from $\wedge$ and $\vee$. This looks somewhat strange, since we showed semantically  that the converse of $S_1$ and $S_2$ are not tautologies. However this cannot be formulated in the  straightforward way, namely   as the {\em schemes} $\varphi|\psi\not\rightarrow \varphi\wedge\psi$ and $ \varphi\vee\psi\not\rightarrow\varphi|\psi$ (the latter are false, e.g. for $\varphi=\psi$). It means that the axiomatic systems $K_i$, for $i=0,1,2,3$ introduced above are {\em interpretable} in the standard propositional logic PL, through the obvious interpretations $I_\wedge$ and $I_\vee$  that interpret $|$ as $\wedge$ or $\vee$, respectively. These are defined inductively in the obvious way for standard connectives, while for $|$ we have   $(\varphi|\psi)^{I_\wedge}=\varphi\wedge\psi$ and $(\varphi|\psi)^{I_\vee}=\varphi\vee\psi$.   Then clearly  for any $\varphi$  and for $I$ being some of these interpretations, $$\vdash_{K_i}\varphi \ \Rightarrow \ \vdash\varphi^I,$$
that is, for every  $Dec$-tautology $\varphi$  (to consider the strongest system $Dec$ of choice functions), $\varphi^I$ is a classical tautology.
However both of the aforementioned  interpretations are not ``faithful'', which means that  the converse of the above implication is not true. For example for classical sentences $\alpha$, $\beta$, $(\alpha|\beta)^{I_\wedge}=\alpha\wedge\beta$,  hence
$(\alpha|\beta\rightarrow \alpha)^{I_\wedge}=\alpha\wedge\beta\rightarrow \alpha$. Then  $\alpha\wedge\beta\rightarrow \alpha$ is a classical tautology
while $\alpha|\beta\rightarrow \alpha$ is not a $K_i$-tautology.

The question is if there exist any further axioms, appropriate for some finer class $X\subset Dec$, which can distinguish $|$ from $\wedge$ and/or $\vee$.  The answer is yes. For example a further condition  that can be imposed to $\neg$-decreasing orderings is one that concerns the  position of the special  classes $\bot$ and $\top$ in this ordering. For example we may require that our decreasing orderings $<$ satisfy the condition $\top<\bot$. Let $Dec_{\top<\bot}$ denote the class of these total orderings of $Sen(L)$. It is rather straightforward that the additional axiom needed to characterize  $Dec_{\top<\bot}$  is

($S_6$) \ $\bot|\top$.

\noindent Thus  $\bot|\top$ is a $Dec_\bot$-tautology while  $(\bot|\top)^{I_\wedge}=\bot\wedge\top$ is not a classical tautology, which means that the logic corresponding to $Dec_{\top<\bot}$ is not interpretable in PL through $I_\wedge$. We might also require that $\top$ is the least element of $<$. If this class of orderings is denoted by $Dec_\top$, the  additional scheme needed to characterize $Dec_\top$ is

($S_6'$) \ $\alpha|\top$.

\noindent Again $Dec_{\top<\bot}$ is not interpretable in PL through $I_\wedge$. Similarly if $Dec_{\bot<\top}$, $Dec_\bot$ denote the classes of decreasing orderings $<$ such that  $\bot<\top$ and $\bot<\alpha$, for every $\alpha$, respectively, then the needed corresponding axioms
are

($S_7$) \ $\neg(\bot|\top)$

\noindent and

($S_7'$) \ $\neg(\alpha|\bot)$,

\noindent respectively.  These logics are not interpretable in PL through $I_\vee$.

\section{Future work}
Our future work will focus on two goals. The first goal is to develop some alternative semantics for superposition logics. We have already  found a second semantics based again on choice functions, but this time the choice applies not to pairs of sentences but to pairs of elements of a Boolean algebra ${\cal B}$, in which the standard sentences  of PL take truth values. We can refer to this  as ``Boolean-value choice semantics'' (BCS) to distinguish it from the ``sentence choice semantics'' (SCS)  used in the present paper.

The second goal is to extend  PLS to {\em first-order superposition logic} (abbreviated FOLS). Such an extension might help us to pass from superposition of sentences/formulas to {\em superposition of objects}. Given two objects (constants) $a$ and $b$, let us consider the formula (in one free variable) $(x=a)|(x=b)$. If our logic can prove that there exists a {\em unique} object $c$ satisfying this formula,  then we can set $c=a\!\uparrow\!\! b$ and say that  $c$ is the  {\em superposition of $a$ and $b$}.  Thus in order for such an operation to be defined for all objects $x,y$ the sentence $(\forall x,y)(\exists!z)((z=x)|(z=y))$ must be a tautology. So the question is whether there is a suitable formalization of FOLS in which this sentence can be stated and  proved.

\vskip 0.2in

{\bf Acknowledgements} Many thanks to two anonymous referees for valuable remarks and suggestions.

\end{document}